\RequirePackage{etoolbox}
\csdef{input@path}{{style/}{graphics/}}
\documentclass[aap,MSNbibl,dvips]{arximspdf}
\makeatletter
   \@ifpackageloaded{graphicx}{}{\usepackage{graphicx}}
\makeatother

\doi{10.1214/14-AAP1023} 
\volume{25}
\issue{3}
\pubyear{2015}
\firstpage{1279}
\lastpage{1324}
\docsubty{FLA}

\makeatletter
\newcommand{\rrvert}{\vert}
\newcommand{\rrVert}{\Vert}
\newcommand{\llvert}{\vert}
\newcommand{\llVert}{\Vert}
\renewcommand{\mid}{|}
\newtheorem{theorem}{Theorem}[section]
\newtheorem{lemma}[theorem]{Lemma}
\newcommand{\nats}{\mathbb{N}}
\newcommand{\ints}{\mathbb{Z}}
\newcommand{\ind}{\mathbb{I}}
\newcommand{\pr}{\mathbb P}
\newcommand{\po}{\operatorname{Po}}
\makeatother

\begin{document}
\begin{frontmatter}

\title{Balanced routing of random calls}
\runtitle{Balanced routing of random calls}

\begin{aug}
\author[A]{\fnms{Malwina J.}~\snm{Luczak}\corref{}\ead[label=e1]{m.luczak@qmul.ac.uk}\thanksref{T1}}
\and
\author[B]{\fnms{Colin}~\snm{McDiarmid}\ead[label=e2]{cmcd@stats.ox.ac.uk}}
\runauthor{M.~J. Luczak and C. McDiarmid}
\affiliation{Queen Mary University of London and University of Oxford}
\address[A]{School of Mathematical Sciences\\
Queen Mary University of London\\
Mile End Road\\
London E1 4NS\\
United Kingdom\\
\printead{e1}} 
\address[B]{Department of Statistics\\
University of Oxford\\
1 South Parks Road\\
Oxford OX1 3TG\\
United Kingdom\\
\printead{e2}}
\end{aug}
\thankstext{T1}{Supported by an EPSRC Leadership Fellowship, Grant
reference EP/J004022/2.}

\received{\smonth{3} \syear{2011}}
\revised{\smonth{12} \syear{2012}}

%
\begin{abstract}
We consider an online network routing problem in continuous time, where
calls have Poisson arrivals and exponential durations. The first-fit
dynamic alternative routing
algorithm sequentially selects up to $d$ random two-link routes between the
two endpoints of a call, via an intermediate node, and assigns the call
to the first route with spare capacity on each link, if there is such a
route. The balanced dynamic alternative routing algorithm
simultaneously selects $d$ random two-link routes, and the call is
accepted on a route minimising the maximum of the loads on its two
links, provided neither of these two links is saturated.

We determine the capacities needed for these algorithms to route calls
successfully and find that the balanced algorithm requires a much
smaller capacity.
In order to handle such interacting random processes on networks, we
develop appropriate tools such as lemmas on biased random walks.
\end{abstract}

%
\end{frontmatter}


\section{Introduction}\label{secintro}

Modern telecommunication systems operate at high\break bandwidth and throughput
and require quick path selection algorithms in order to fully utilise
network resources while minimising routing cost. In many settings, each pair
of nodes have dedicated capacity for communication between them,
designed to
meet demand. When all the capacity is in use in times of congestion, common
routing strategies will attempt to find an alternative route via one or
more intermediate nodes. Usually, an admission protocol checks a small number
of alternatives, and rejects the incoming call if none is available. Examples
of such protocols include AT\&T's Dynamic Nonhierarchical Routing
algorithm~\cite{acm} and the Dynamic Alternative Routing (DAR)
algorithm~\cite{gkk}; see also~\cite{ch,ghk,h,hl,k91}.

Dynamic routing in communication networks belongs to a class of online
load-balancing problems, where tasks are to be
assigned to one or more links (servers), and communication requests (customers)
may only be assigned to specific paths (subsets of servers), depending on
their properties and/or network topology. Research in this area has witnessed
rapid developments, with many papers demonstrating the advantage of balanced
allocations, as in the ``power of two choices''
phenomenon~\cite{bcsv99,klm96,lmc03a,lmc03b,lmcu03,lu99,m96,mrs01}.

This paper is concerned with an online routing problem in continuous time,
where calls have Poisson arrivals and exponential durations
(and so in particular calls end, in contrast to many earlier models).
Load-balancing and alternative routing strategies are deployed to assign
bandwidth to arriving calls, under constraints imposed by network topology.
First, in order to set the scene, let us recall a related online routing
problem in discrete time from~\cite{lmcu03}, where calls do not end.

\textit{An earlier discrete time model.}
There is a set $V=\{1, \ldots, n\}$ of $n$ nodes,
each pair of which may wish to communicate. A \textit{call} is an
unordered pair $\{u,v\}$ of distinct nodes, that is an edge of the
complete graph $K_n$ on $V$. For each of the
${n \choose2}$
unordered pairs $\{u,v\}$ of distinct nodes, there is a \textit{direct
link}, also denoted by $\{u,v\}$, with capacity $D_1=D_1(n)$. The
direct link is used to route a call
as long as it has available capacity. There are
also two \emph{indirect links}, denoted by $uv$ and $vu$, each with
capacity $D_2=D_2(n)$. The indirect link $uv$ may be used when for
some $w$ a call $\{u,w\}$ finds its direct link saturated, and we
seek an alternative route via node $v$. Similarly $vu$ may be used
for alternative routes for calls $\{v,w\}$ via $u$.

We are given a sequence of $M$ calls one at a time. For
each call in turn, we must choose a route (either a direct link or an
alternative two-link route via an intermediate node) if this is possible,
before seeing later calls. These routes cannot be changed later, and
calls do
not end. The aim is to minimise the number of calls that fail to be routed
successfully.

The calls are independent random variables
$Z_1,Z_2, \ldots, Z_M $, where each $Z_j$ is uniformly distributed
over the
edges $e \in E(K_n)$, the edge set of $K_n$.
Let $d$ be a (fixed) positive integer.
A \emph{general dynamic alternative routing algorithm} GDAR operates as follows.
For each call $e = \{u,v\}$ in turn, the call is routed on the direct link
if possible, and otherwise nodes $w_1, \ldots, w_d$ are selected
uniformly at random with replacement from $V \setminus\{u,v\}$,
and the call is routed via one of these nodes if possible, along
the two corresponding indirect links. The \emph{first-fit dynamic
alternative routing algorithm} FDAR is the version when we always
choose the first possible alternative route,
if there is one. The \emph{balanced dynamic alternative
routing algorithm} BDAR
is the version when we choose an alternative
route which minimises the larger of the current loads on its two
indirect links, if possible.
Calls that do not find an available route are lost.

Results for this model were first obtained in~\cite{l00,lu99},
and later strengthened and extended in~\cite{lmcu03}. Consider the
case where
$M \sim c{n \choose2}$ for a constant $c>0$. It is known that with the
algorithm FDAR we need
both link capacities $D_1,D_2$ of order $\sqrt{{ \ln n} \over\ln\ln n}$
to ensure that \emph{asymptotically almost surely} (a.a.s.), that is, ``with
probability $\rightarrow1$ as $n \rightarrow\infty$'', all $M$ calls are
routed successfully.
The balanced
method BDAR succeeds with much smaller capacities. Specifically,
there is a tight threshold value close to $\ln\ln n/\ln d$ for $D_2$ to
guarantee that a.a.s. no call fails (and
the precise value of $D_1$ is unimportant;
see Theorems 1.3~and~7.1 in~\cite{lmcu03}, where in the latter $D_1=0$).


\subsection{Our model} \label{subsecmodel}

Here we consider a related continuous-time network model,
with the desirable additional feature that calls end.
Of course this gives a much better model for calls, but it leads to
harder analysis,
since, for example, we now need to handle biased random walks with
negative as well as positive increments.

Calls arrive in a Poisson process with rate
$\lambda{n \choose2}$, where $\lambda$ is a positive constant.
The calls are i.i.d. random variables $Z_{1},Z_{2}, \ldots,$ where
$Z_{j}$ is the $j$th call to arrive and is uniform over the
edges of $K_n$ for each $j$; also let $T_j$ be the arrival time of call $Z_{j}$.
For each edge $\{u,v\}$ there are two links, $uv$ and $vu$, both with
capacity $D=D(n)< \infty$. Since in~\cite{lmcu03} the use of direct
links was found to have only a minor effect on the total capacity
requirements for efficient communication, here we do not use direct
links but instead demand that each call be routed along a path
consisting of a pair of indirect links.
This yields a cleaner model which captures the interesting behaviour,
and for which we can give
a rigorous analysis without (we hope!), making the paper
too long for the gentle reader.

If a call is for $\{u,v\}$, then we pick $d$ possible intermediate nodes
uniformly at random with replacement, as in the GDAR algorithm. The FDAR
algorithm chooses the first possible alternative route, if
there is one.
The BDAR algorithm chooses an alternative route minimising the larger
of the current loads
on its two links, if possible (ties are broken arbitrarily).
Call durations are unit mean
exponential random variables, independent of one another and of the
arrivals and choices processes. When a call terminates, both
busy links are freed. Calls that do not find an available route are lost.

For each edge $e=\{u,v\} \in E(K_n)$ and node $w \in V \setminus e$,
let $X_t (e,w)$
denote the number of calls in progress at time $t$ which are
routed along the path
consisting of links $uw$ and $vw$, that is, calls between the end nodes
$u$ and
$v$ of $e$ routed via $w$.
We call $X_t = (X_t (e,w)\dvtx  e \in E, w \in V \setminus e)$ the \emph{load
vector} at time $t$
and let ${\mathcal X} = (\ints^+)^{n(n-1)(n-2)/2}$ denote the set of
all possible load vectors.
The process $X=(X_t)_{t \ge0}$ of load vectors is
a continuous-time
jump Markov chain with state space~${\mathcal X}$,
defined on some probability space $(\Omega, \mathcal F, \pr)$.
By standard results, there exists a
unique stationary distribution $\pi$, and, whatever the
distribution of the starting state $X_0$, the distribution of the
load vector $X_t$ at time $t$ converges to $\pi$ as
$t \rightarrow\infty$.

We put a natural partial order on ${\mathcal X}$: given two vectors $x,
\tilde{x} \in{\mathcal X}$, we say that $x \le\tilde{x}$ if $x(e,w)
\le\tilde{x}(e,w)$ for each $e \in E(K_n)$, $w \in V\setminus e$.
Given \mbox{${\mathcal X}$-}valued random variables $Z$ and $\tilde{Z}$, we
say that $\tilde Z$ stochastically dominates
$Z$ if $\pr(Z \ge z) \le \pr(\tilde{Z} \ge z)$ for all $z$. If this
is the case, then we also say that the distribution $F_Z$ of $Z$ is
stochastically dominated by the distribution $F_{\tilde Z}$ of $\tilde
Z$. We note that $\tilde Z$ stochastically\vspace*{1pt} dominates $Z$ if and only if
there exists a coupling of $Z$ and $\tilde{Z}$ such that $Z \le\tilde
{Z}$ with probability 1.

Our main interest is in the \emph{blocking probability}, that is,
the probability that a new
call fails to find an available route and is thus lost.
As in the discrete version analysed in~\cite{lmcu03}, or in
the models analysed in~\cite{lmc03a} and~\cite{lmc03b}
(see also~\cite{abku94,bcsv99,m96,mrs01}), in our more complicated
continuous-time network model
we observe the ``power of two choices'' phenomenon; that is, with the
BDAR algorithm for $d \ge2$ the capacity required to ensure that
most calls are routed successfully is much smaller
than with the FDAR algorithm.
(When $d=1$ FDAR and BDAR reduce to the same algorithm.)
Let us now state our main results, which contain precise statements of
this maxim.

Throughout the paper, we use the asymptotic $O(\cdot)$, $\Omega(\cdot)$ and $o(\cdot)$ notation in a usual way.
Thus for nonnegative functions $f(n)$ and $g(n)$ defined on $\nats$,
we write $f(n) = O(g(n))$
if there exists a constant $C$ such that $f(n) \le C g(n)$ for all
sufficiently large $n$,
$f(n) = \Omega(g(n))$ if $g(n) = O(f(n))$, and $f(n) = o(g(n))$ if
$f(n)/g(n) \to0$ as $n \to\infty$.
%


\subsection{Our results} \label{subsecresults}

Theorem~\ref{thmfdar1} below shows that, when the FDAR algorithm is used,
capacity $D(n)$ of order ${{\ln n} \over{\ln\ln n}}$ is needed in order
to ensure that no call is lost in a time interval of length polynomial
in $n$.
The set-up is as follows.

The arrival rate per edge is fixed as $\lambda>0$, and $d$ is a fixed
positive integer.
Let $\alpha>0$, and let each link have capacity $D = D(n) \sim\alpha
{{\ln n} \over{\ln\ln n}}$
as $n \to\infty$. We may need a ``burn-in'' period $t_0$:
for each $n$, if the distribution of the initial state $X_0$ is
stochastically dominated by the stationary distribution $\pi$,
then let $t_0=0$, and otherwise let $t_0=t_0(n)= 5 \ln n$.
Now we consider any $t_1 \geq t_0$ and $K>0$, and time intervals
$[t_1,t_1+n^K]$.

Let us say that $\alpha$ is $K$-\emph{good} if, whatever version of
GDAR we use, for each $t_1 \geq t_0$,
the mean number of calls lost during the interval $[t_1, t_1 + n^K]$ is $o(1)$;
and $\alpha$ is $K$-\emph{bad} if, when we use FDAR, for each $t_1
\geq0$,
the mean number of calls lost during the interval $[t_1, t_1 + n^K]$ is
$n^{\Omega(1)}$.
(Observe that $\alpha$ cannot be both $K$-good and $K$-bad.)

\begin{figure}

\includegraphics{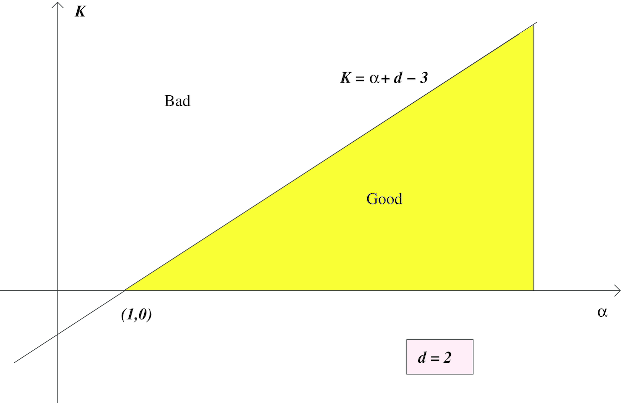}

\caption{When $\alpha$ is $K$-good: case $d \le2$.} \label{fig2}
\end{figure}

\begin{figure}

\includegraphics{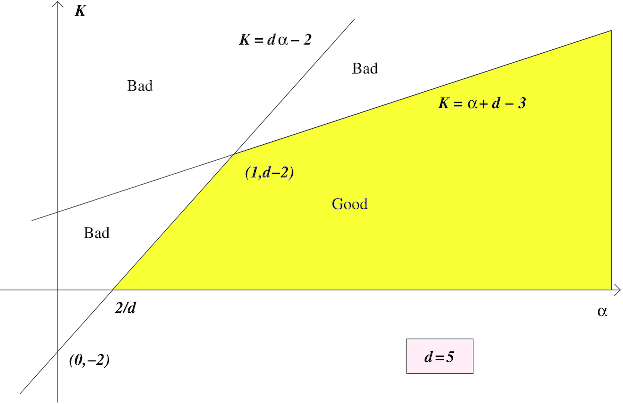}

\caption{When $\alpha$ is $K$-good: case $d > 2$.} \label{fig1}
\end{figure}

The first theorem below shows that $\alpha=2/d$ is a critical value
(which does not depend on $\lambda$). In particular, if $\alpha>2/d$, then
$\alpha$ is $K$-good for some $K>0$.
The second theorem concerns $\alpha$ above this threshold and
describes the pairs
$\alpha,K$ where $\alpha$ is $K$-good or $K$-bad.
The behaviour is simple when $d$ is 1 or 2, and more interesting for $d
\geq3$;
see Figures~\ref{fig2} and~\ref{fig1}.

%

%

%
\begin{theorem} \label{thmfdar1}
If $\alpha> 2/d$, then $\alpha$ is $K$-good for some $K>0$, and if
$\alpha\leq2/d$, then $\alpha$ is $K$-bad for each $K>0$.
\end{theorem}

%
\begin{theorem} \label{thmfdar2}
Let $\alpha> 2/d$, and let $K>0$.
\begin{longlist}[(a)]
\item[(a)]
If $2/d< \alpha\leq1$ (and so $d \geq3$), then $\alpha$ is $K$-good
for $d\alpha-K > 2$, and $\alpha$ is $K$-bad for $d\alpha-K < 2$.
\item[(b)]
If $\alpha\geq1$ (as must be the case when $d$ is 1 or 2), then
$\alpha$ is $K$-good
for $\alpha-K > 3-d$, and $\alpha$ is $K$-bad for $\alpha-K < 3-d$.
\end{longlist}
\end{theorem}

As foreshadowed above, the next result
shows that the BDAR algorithm requires significantly smaller capacities.
Note that the expected number of calls arriving in a time interval of
length $n^K$
is ${\sim}(\lambda/2) n^{K+2}$.
%

\begin{theorem} \label{thmbdar}
Let $\lambda>0$ be fixed, and let $d \geq2$ be a fixed integer.
Let $K >0$ be a constant.
Then there exist constants $c=c(\lambda,d,K)>0$ and $\kappa= \kappa
(\lambda,d)$
such that the following holds:
\begin{longlist}[(s)]
\item[(a)]
Suppose that $D(n) \ge{\ln\ln n}/ {\ln d} +c$, and we use the BDAR algorithm.
Given $n$ and a distribution for $X_0$, let $t_0=0$ if this initial distribution
is stochastically dominated by the equilibrium distribution, and let
$t_0 =\kappa\ln n$ otherwise.
Then the expected number of failing calls during the interval $[t_1,
t_1 + n^K]$ is $o(1)$
for each $t_1 \geq t_0$.
\item[(b)]
If\vspace*{1pt} $D(n)\le{\ln\ln n}/{\ln d} -c$ and we use any GDAR algorithm,
then a.a.s. at least $n^{K+2- o(1)}$ calls are lost during
$[t_1, t_1 + n^K]$
for each $t_1 \geq0$.
\end{longlist}
\end{theorem}

Some parts of our proofs are built on our earlier work on balls and
bins in continuous time~\cite{lmc03a}.
Indeed that earlier paper arose from the need of the authors to sort
out simpler ``network-free''
load-balancing results so as to be ready to
tackle the additional complications in network routing problems, where
a call occupies two adjacent links.

We mention that a process similar to the one defined above, but
sometimes also with direct links,
was considered in Luczak and Upfal~\cite{lu99}
and then in Anagnostopoulos, Kontoyiannis and Upfal~\cite{aku}.
The earlier of these works obtained, heuristically, some preliminary results.
These indicate that link capacity of order $\ln\ln n/\ln d$ is
sufficient to ensure that with the BDAR algorithm,
in equilibrium, a new call is accepted with high probability,
and capacity of order $\ln\ln(t_0n)/\ln d$ is sufficient to ensure
that all calls arriving during an
interval of length $t_0$ are accepted with high probability.
There is also a short explanation of why link capacity of order $\Omega
(\sqrt{\ln(t_0n)/\ln\ln(t_0n)})$
is necessary to achieve this with FDAR.

Augmented
versions of these arguments appear in the later paper of
Anagnostopoulos et al.~\cite{aku}.
They find an upper bound of $\ln\ln n/\ln d + o(\ln\ln n/\ln d)$ for
the capacity required by the
BDAR algorithm to ensure that, in equilibrium, an arriving call is
accepted with probability
tending to 1 as $n \to\infty$.
Further, they identify a lower bound of $\Omega(\sqrt{\ln n/\ln\ln n})$
for the capacity needed by the FDAR algorithm to achieve the same effect.

Here we give rigorous proofs of sharp
versions of these bounds and turn them into sharp two-sided results by
supplementing them with a matching
lower bound on the performance of the BDAR algorithm
and a matching upper bound on the performance of the FDAR algorithm.
Further, we do not restrict our attention to the equilibrium
distribution, and we
prove upper and lower bounds on the performance of these algorithms
over long time intervals.
Accordingly, our proofs are considerably more involved and subtle than
the arguments
put forward
in~\cite{aku}. In comparison with that paper,
our lower bounds for the FDAR algorithm are of the order $\ln n/\ln\ln
n$, not
$\sqrt{\ln n/\ln\ln n}$: this is due
to the fact that we never allow direct routing between pairs of nodes,
whereas in~\cite{aku} direct routing is allowed for FDAR (though not
for BDAR).


\subsection{Some notation} \label{subsecnotation}

Here we give some further definitions and notation which we shall need shortly.
The subscript $t$ always refers to time.
Given an edge $e = \{u,v\} \in E(K_n)$, let $X_t(e) = \sum_{w \notin
e} X_t(e,w)$
denote the number of calls between $u$ and $v$ in progress at time $t$.
Also, given distinct nodes $v$ and~$u$, let
$X_t (vu)= \sum_{w \neq u,v} X_t (\{v,w\},u)$, which is the load of link~$vu$
at time $t$.
Given a node~$v \in V$, let $X_t(v)=\sum_{u \neq v} X_t (vu)$,
which is the number of calls with one end~$v$ at time~$t$.
Thus $\llVert X_t \rrVert_1:={1 \over2} \sum_{v \in V}
X_t(v)$ is the total
number of calls at time~$t$.

We say that a link is \emph{saturated} (or full) if it has load equal to
its capacity~$D$.
Given a node~$v$, we let ${\mathcal S}_{t}(\mbox{at } v) $ denote the
set of saturated
links $vw$ (for $w \neq v$) for calls
at~$v$ at time~$t$
and let
$S_{t}(\mbox{at } v) = \llvert{\mathcal S}_{t}(\mbox{at } v) \rrvert
$, which is
the number of
saturated links $vw$.
Similarly, given a node~$w$, we let ${\mathcal S}_{t}(\mbox{via } w) $
denote the set of
saturated links $vw$
for calls at some node $v \neq w$ at time~$t$ and let $S_{t}^{}(\mbox
{via } w) =\llvert{\mathcal S}_{t}(\mbox{via } w) \rrvert$.


\subsection{Overview of the proofs} \label{overview}

The rest of the paper is organised as follows.
Section~\ref{secprelim} contains some preliminary lemmas that will be
needed later in our proofs.
After recalling some probability inequalities, there are results
concerning biased random walks, transferring probability bounds from
points to intervals,
comparing jump Markov chains to independent birth-death processes,
and a nonasymptotic version of the PASTA principle.
Section~\ref{sectotal} is where we introduce the ``network'' dependencies.
Lemma~\ref{lemblocking} is a key result on the probability of a call
failing conditionally
on the history up to then, and also we establish
inequalities for the total number of calls and the number of saturated links.

In Section~\ref{secproof1} we prove Theorems~\ref{thmfdar1}~and~\ref{thmfdar2}, which
describe when the pair $(\alpha,K)$ is good or bad.
The approximate picture is as follows.
The number $S_{t}(\mbox{at } v) $ of saturated links at a node $v$ has
expected value $n^{1-\alpha+ o(1)}$, and
the probability $p$ that a call with one end $v$ fails is roughly\vspace*{1.5pt}
$\mathbb E[(S_{t}(\mbox{at } v) /n)^d]$.
If $0<\alpha<1$, then $S_{t}(\mbox{at } v) $ is concentrated and
$\mathbb E[S_{t}(\mbox{at } v) ^d] = n^{(1-\alpha)d + o(1)}$ and $p$
is $n^{-\alpha d + o(1)}$.
The\vspace*{1.5pt} expected number of arrivals in an interval of length $n^K$ is about
$n^{K+2}$,
and $n^{K+2} p = n^{K+2 - \alpha d +o(1)}$, so $\alpha$ is $K$-good
when $K+2 - \alpha d <0$,
and $\alpha$ is \mbox{$K$-}bad when $K+2 - \alpha d >0$.
When\vspace*{1pt} $\alpha\geq1$, then $\mathbb E[S_{t}(\mbox{at } v) ^d] \sim
\mathbb E[S_{t}(\mbox{at } v) ]$ and
$p = n^{1-\alpha-d +o(1)}$, and again we see when $\alpha$ is
$K$-good by looking at $n^{K+2} p$.

To show goodness in these theorems we need to show that, throughout a
time interval,
there are not too many saturated links in the network.
To show badness we need a lower bound on the number of saturated links
at a vertex,
and for this we need also to upper bound the number of saturated links,
so that an arriving call wishing
to use the link is not too often blocked because the ``partner'' link
of the pair is saturated.

In Section~\ref{secproof2} we prove Theorem~\ref{thmbdar} on the
balanced routing algorithm BDAR.
We are not able to use the neat approach used in~\cite{lmc03a} for
balls in bins
(based on rapid mixing, concentration and simple explicit balance
equations in equilibrium)
in the more complicated network model. Instead the proof is based on
the ``layered induction'' approach,
used, for example, in~\cite{abku94,abku99},
though now with additional hurdles.

For the upper bound, the key step is to show that if for each node $v$
the number of arcs at $v$
with ``weighted load'' at least $h$ is at most $\alpha$ throughout an
interval $[t, t_0]$,
then with high probability for each node $v$ the number of arcs at $v$
with weighted load at least
$h+1$ is at most $\alpha' \ll \alpha$ throughout a slightly smaller
interval $[t', t_0]$.
We thus deduce that with high probability no link is ever saturated in
the relevant interval of length $n^K$,
and so no call is lost.
For the lower bound, we use a similar approach to show that with high
probability
for each node $v$, at least $n^{1-\varepsilon}$ links $vw$ incident on
$v$ are saturated
throughout the interval of length $n^K$, and hence with high
probability at least
$ n^{K +2-\varepsilon d-o(1)}$ calls arriving during the interval are lost.

Finally we make some concluding remarks in Section~\ref{secconc}.


\section{Preliminary results}\label{secprelim}

In this section we give some basic results which will be used in our proofs.
Topics covered include some general probability inequalities and random
walks ``with drift''.
The reader may wish to skim this section and refer back to it as required.

\subsection{Inequalities}\label{subineq}

If the random variable $Y$ has the Poisson distribution with mean $\mu
>0$, we write
$Y \sim\po(\mu)$, and for nonnegative integers $k$, we write
%
%
\begin{equation}
\label{eqnpobound2} p_k(\mu)= \pr(Y \ge k) = e^{-\mu} \sum
_{j \ge k}{{\mu^{j}}\over{j!}}
\end{equation}
and note that
%
%
\begin{equation}
\label{eqnpobound1} p_k(\mu) \le\mu^k/k! \le(e
\mu/k)^k.
\end{equation}
When $\mu>0$ is a constant and $D=D(n)$ is an integer with $D \sim
\alpha\ln n / \ln\ln n$, we have
%
%
\begin{equation}
\label{eqnDpo} p_D(\mu) = n^{-\alpha+o(1)}.
\end{equation}
The following are a pair of standard concentration inequalities for a
binomial or
Poisson random variable $Y$ with mean $\mu$:
%
%
\begin{equation}
\label{binpo1} \pr(Y-\mu\ge\varepsilon\mu) \le\exp\bigl(-\tfrac13
\varepsilon^2 \mu\bigr)
\end{equation}
and
%
%
\begin{equation}
\label{binpo2} \pr(Y -\mu\le-\varepsilon\mu) \le\exp\bigl(-\tfrac12
\varepsilon^2 \mu\bigr)
\end{equation}
for $0 \le\varepsilon\le1$; see, for example, Theorem~2.3(c) and inequality~(2.8) in~\cite{cmcd98}.

We shall use the following version of Talagrand's inequality; see, for example,
Theorem~4.3 in~\cite{cmcd98}. (In the notation in~\cite{cmcd98}, the function
$h$ below is a \mbox{$(c^2 r)$-}configuration function.)

%
\begin{lemma} \label{lemtalagrand}
Let ${\mathbf Y} = (Y_1,Y_2,\ldots)$ be a finite family of independent random
variables, where the random variable $Y_j$ takes values in a set
${\mathcal Y}_j$. Let ${\mathcal Y} = \prod_j {\mathcal Y}_j$.

Let $c$ and $r$ be positive constants, and suppose that the
nonnegative real-valued measurable
function $h$ on ${\mathcal Y}$ satisfies the following two
conditions for each ${\mathbf y }\in{\mathcal Y}$:
\begin{itemize}
\item
Changing the value of a co-ordinate $y_j$ can change the value of
$h({\mathbf y})$ by at most~$c$.
\item
If $h({\mathbf y}) = s$, then there is a set of at most $rs$ co-ordinates
such that $h({\mathbf y'}) \geq s$ for any ${\mathbf y'} \in{\mathcal Y}$ which
agrees with ${\mathbf y}$ on these co-ordinates.
\end{itemize}
Let $m$ be a median of the random variable $Z=h({\mathbf Y})$.
Then for each $x \geq0$,
%
%
\begin{equation}
\label{above} \pr(Z \geq m+x) \leq2 \exp\biggl(-\frac{x^2}{4 c^2
r(m+x)} \biggr)
\end{equation}
and%
\begin{equation}
\label{below} \pr(Z \leq m-x) \leq2 \exp\biggl(-\frac{x^2}{4 c^2 rm}
\biggr).
\end{equation}
\end{lemma}

\subsection{Random walks and birth-and-death processes}\label{subrw}

We start with a lemma from~\cite{lmcu03}, which will be used in
Sections~\ref{secupper} and~\ref{seclower}.
For $n \in\nats$ and $0 \le p \le1$, let $B(n,p)$ denote a binomial
random variable with parameters $n$ and $p$.

%
\begin{lemma}[(Lemma~2.3 in~\cite{lmcu03})]
\label{lembin}
Let $\mathcal F_0 \subseteq\mathcal F_1 \subseteq\cdots$ be a filtration;
let $Y_1,Y_2, \ldots$ be indicator random variables such that each
$Y_i$ is $\mathcal F_i$-measurable;
and let $E_0,E_1, \ldots$ be events where $E_i \in\mathcal F_i$ for
each $i=0,1,\ldots.$
For each $t \in\nats$, let $R_t = \sum_{i=1}^t Y_i$. Let $0 \le p
\le1$, and let $k$ be a positive integer.
\begin{longlist}[(a)]
\item[(a)] Assume that for each $i=1,2, \ldots$
\[
\pr(Y_i = 1 \mid\mathcal F_{i-1}) \le p\qquad\mbox{on }
E_{i-1} \cap\{R_{i-1} < k\}.
\]
Then for each $t \in\nats$
\[
\pr\Biggl(\{R_t \ge k\} \cap\Biggl( \bigcap
_{i=0}^{t-1} E_i \Biggr) \Biggr) \le\pr
\bigl(B(n,p) \ge k\bigr).
\]
\item[(b)] Assume that for each $i=1,2, \ldots$
\[
\pr(Y_i = 1 \mid\mathcal F_{i-1}) \ge p\qquad\mbox{on }
E_{i-1} \cap\{R_{i-1} < k\}.
\]
Then for each $t \in\nats$
\[
\pr\Biggl(\{R_t < k\} \cap\Biggl( \bigcap
_{i=0}^{t-1} E_i \Biggr) \Biggr) \le\pr
\bigl(B(n,p) < k\bigr).
\]
\end{longlist}
\end{lemma}

The next lemma concerns hitting times of a generalised random walk
with a ``downward drift''. It is the ``reverse'' of Lemma~7.2 in~\cite{lmc03a},
and can be deduced easily from that result by replacing the $Y_i$ with $-Y_i$;
we omit the details.
It will be used in the proofs of Lemma~\ref{lemhittime1} and
Theorem~\ref{thmbdar}(b).

%
\begin{lemma} \label{lemhittime}
Let $\mathcal F_0 \subseteq\mathcal F_1 \subseteq\cdots$ be a filtration;
let $Y_1,Y_2, \ldots$ be random variables taking values
in $\{-1,0,1\}$ such that each $Y_i$ is $\mathcal F_i$-measurable;
and let $E_0, E_1, \ldots,$ be events where each $E_i \in\mathcal F_i$.
For each $t \in\nats$, let $R_t=R_0 + \sum_{i=1}^t Y_i$.
Let $0 \le p \le1/3$, let $r_0$ and $r_1$ be integers such that $r_1< r_0$,
and let $m$ be an integer such that $pm \geq2(r_0-r_1)$.
Assume that for each $i=1,\ldots,m$,
\[
\pr(Y_i =1\mid\mathcal F_{i-1}) \le p \qquad\mbox{on }
E_{i-1} \cap(R_{i-1} > r_1)
\]
and
\[
\pr(Y_i=-1 \mid\mathcal F_{i-1}) \ge2p \qquad\mbox{on }
E_{i-1} \cap(R_{i-1} > r_1).
\]
Then
\[
\pr\Biggl( \Biggl( \bigcap_{t=1}^{m}
\{R_t > r_1\} \Biggr) \cap\Biggl(\bigcap
_{i=0}^{m-1} E_i\Biggr) \bigg|
R_0=r_0 \Biggr) \le\exp\biggl( - \frac{pm}{28}
\biggr).
\]
\end{lemma}

We will use the last lemma to show that, for a type of discrete-time
``immigra\-tion-death'' process
satisfying suitable conditions, it is unlikely that the ``population''
$R_t$ stays above
the level $r$ throughout a long period.
The following lemma will be used in the proof of Theorem~\ref{thmbdar}(a).

\begin{lemma} \label{lemhittime1}
Let $\mathcal F_0 \subseteq\mathcal F_1 \subseteq\cdots,$ be a
filtration; let
$Y_1,Y_2, \ldots$ be random variables taking values in $\{-1,0,1\}$ such
that each $Y_i$ is $\mathcal F_i$-measurable; and let $E_0, E_1, \ldots
$ be
events where each $E_i \in\mathcal F_i$.
Let $a,b>0$ be constants, and let $\tilde{r}$~and~$r$ be integers with
$2a/b \le r \leq\tilde{r}-1$.

Let $R_0 = \tilde{r}$, and let $R_t=R_0 + \sum_{i=1}^t Y_i$. Assume that
for each $i=1,2,\ldots$
\[
\pr(Y_i =1\mid\mathcal F_{i-1}) \le a\qquad\mbox{on }
E_{i-1} \cap(R_{i-1} >r)
\]
and
\[
\pr(Y_i=-1\mid\mathcal F_{i-1})\ge by\qquad\mbox{on }
E_{i-1} \cap(R_{i-1}=y)
\]
for each $y=r+1,\ldots,\tilde{r}$, and
\[
\pr(Y_i=-1\mid\mathcal F_{i-1})\ge b\tilde{r}\qquad\mbox{on
} E_{i-1} \cap(R_{i-1}>\tilde{r}).
\]
Let $m'= \lceil\frac{4}{b} \rceil\lceil\log_2 \frac{\tilde
{r}}{r} \rceil$,
and let $E$ be the event $\bigcap_{i=1}^{m'} E_i$.
Then
%
%
\begin{equation}
\label{eqhittime} \pr\Biggl(\Biggl( \bigcap_{t=1}^{m'}
\{R_t > r\} \Biggr) \cap E\Biggr) \le2\exp\biggl(-
{r\over14} \biggr).
\end{equation}
\end{lemma}

\begin{pf}
Let $k = \lceil\log_2 \frac{\tilde{r}}{r} \rceil-1$, so that
$2^k r < \tilde{r} \leq2^{k+1} r$.
Let $T_0, T_1,\ldots,T_k$ be the hitting times to cross the $k+1$
intervals from $\tilde{r}$ down to $2^kr$, from $2^kr$ down to
$2^{k-1}r$, and
so on, ending with the interval from $2r$ down to $r$. Thus
\[
T_0 = \min\bigl\{ t \geq0\dvtx  R_{t} = 2^{k}r
\bigr\},
\]
and for $j=1,\ldots,k$,
\[
T_j = \min\bigl\{ t > T_{j-1}\dvtx  R_{t} =
2^{k-j}r \bigr\}.
\]

Consider\vspace*{2pt} $j \in\{0,\ldots,k\}$.
We want to upper bound the probability that $T_j - T_{j-1}$ is large.
To do this, we may use Lemma~\ref{lemhittime} with $p$ as $p_j = b
2^{k-j-1}r$,
$r_0 = 2^{k-j+1}r$ (except that for $j=0$ we let $r_0=\tilde{r}$),
$r_1 = 2^{k-j}r$ and $m$ as $m_j = \lceil\frac{4}{b} \rceil$.
Note that $p_j m_j \geq2^{k-j+1}r$, which is at least twice the length
of the interval.
(It~may look at first sight that we are ``giving away'' rather a lot on the
``upward'' probability but this makes only a constant factor difference.)
Hence, with $T_{-1} \equiv0$,
\[
\pr\bigl(E \cap\{ T_j-T_{j-1} > m_j\}
\bigr) \leq\exp\biggl( - \frac{p_j
m_j}{28} \biggr) \leq\exp\biggl( -
\frac{2^{k-j}r}{14} \biggr).
\]
But now
\begin{eqnarray*}
\pr\bigl(E \cap\bigl\{ R_t > r\ \forall t \in\bigl\{1,
\ldots,m' \bigr\}\bigr\}\bigr) & \leq& \sum
_{j=0}^{k} \pr\bigl(E \cap(T_j-T_{j-1}
> m_j)\bigr)
\\
& \leq& \sum_{j=0}^{k} \exp\biggl( -
\frac{2^{k-j}r}{14} \biggr)
\\
& \leq& e^{-r/14 }/\bigl(1- e^{-r/14}\bigr).
\end{eqnarray*}
Hence
\[
\pr\bigl(E \cap\bigl\{ R_t > r\ \forall t \in\{1,\ldots,m \}\bigr\}
\bigr) \leq2 e^{-r/14},
\]
by the above if $e^{-r/14} \leq\frac{1}2$ and trivially otherwise.
\end{pf}

The next lemma appears as Lemma 7.3 in~\cite{lmc03a} and shows that
if we
try to cross an interval against the drift we rarely succeed. It will
be used in the proof
of Lemma~\ref{lemupper2} and also in the proof of Lemma~\ref{lemlower2}.

%
\begin{lemma} \label{lemcrossing}
Let $a$ be a positive integer.
Let $p$ and $q$ be reals with $q > p \geq0$ and $p+q \leq1$.
Let $\mathcal F_0 \subseteq\mathcal F_1 \subseteq\mathcal F_2
\subseteq\cdots$ be a filtration;
let $Y_1,Y_2, \ldots$ be random variables taking values in
$\{-1,0,1\}$ such that each $Y_i$ is $\mathcal F_i$-measurable; and
let $E_0, E_1, \ldots$ be events where each $E_i \in\mathcal F_i$.
Let $R_0 = 0$, and let $R_k = \sum_{i=1}^k Y_i$ for $k=1,2,\ldots.$
Assume that for each $i=1,\ldots,m$,
%
\begin{eqnarray}
\pr(Y_i = 1\mid\mathcal F_{i - 1}) \le p\quad\mbox{and}\quad \pr(Y_i = -1 \mid\mathcal F_{i - 1}) \ge q
\nonumber
\\
\eqntext{\mbox{on }E_{i - 1} \cap\{0 \leq R_{i -
1} \leq a - 1\}.}
\end{eqnarray}
Let
\[
T= \inf\bigl\{ k \geq1\dvtx  R_k \in\{ -1,a \} \bigr\}\quad\mbox{and}\quad
E_T = \bigcap_{i=0}^{T}
E_i.
\]
Then
\[
\pr\bigl( \{ R_T =a\} \cap E_T \bigr)
\leq(p/q)^a.
\]
\end{lemma}

We must handle random processes like $X_t(v)$, the number of active
calls with one end $v$ at time $t$,
which can increase only when new calls arrive,
and be able to move from probability bounds
at points of time to bounds
over intervals of time.
We require another lemma, which extends Lemma 2.1 in~\cite{lmc03a}.

Consider the $n$-node case of our network model, where the set of all
load vectors is
${\mathcal X} =(\ints^+)^{n(n-1)(n-2)/2}$. Let us say that a
real-valued function $f$ on ${\mathcal X}$
has \emph{bounded increase at} a node $v$ if whenever $s$ and $t$ are
times with $s<t$,
then $f(x_t)$ is at most $f(x_s)$ plus the total number of arrivals in the
interval $(s,t]$ for $v$;
$f$~has \emph{bounded increase via} a node $v$ if for each $s<t$,
$f(x_t)$ is at most $f(x_s)$ plus twice the total number of arrivals in the
interval $(s,t]$ routed via $v$ as the intermediate node; and for each $s<t$,
$f$ has \emph{strongly bounded increase at} a node $v$ if $f(x_t)$ is at
most $f(x_s)$ plus the maximum number of arrivals for $v$ in the
interval $(s,t]$ which use any given link incident on $v$.
Thus, for example, given $v \in V$, $f(x)= x(v)$ (the total number of
calls with one end $v$ in state $x$) has bounded increase at $v$,
$f(x)=\llvert\{w \in V \setminus\{v\}\dvtx  x(wv) \ge D\}\rrvert$ (number
of saturated
links $wv$, for calls with one end $w$ routed via $v$, in state $x$)
has bounded increase via $v$,
and $f(x)=\max_{w \in V \setminus\{v\}} x(vw)$ (maximum load of a
link $vw$, for calls with one end $v$, in state $x$) has strongly
bounded increase at $v$.

The following elementary lemma will be invoked many times, in the
proofs of various other lemmas,
as well as in the proofs of the three theorems. (Think of the bounds
$g$ as increasing and $h$ as decreasing.)
%

\begin{lemma} \label{lemstretch}
Consider functions $f\dvtx  {\mathcal X} \rightarrow\mathbb{R}$
and $g,h\dvtx  \mathbb{R}\rightarrow\mathbb{R}$.
Let $v$ be a node in $V$, let $t_1 \geq0$ and $\tau>0$, and let $E
\in\mathcal F_{t_1}$.
Suppose that, for all $a \in\mathbb{R}$ and all times $t_1 \le t \le t_1
+ \tau$,
\[
\pr\bigl(E \cap\bigl\{f(X_t) \le a \bigr\} \bigr) \le g(a)\quad\mbox{and}\quad
\pr\bigl(E \cap\bigl\{f(X_t) \ge a \bigr\} \bigr) \le h(a).
\]
Let $\sigma>0$, let $a \in\mathbb{R}$, and let $b \geq0$.
\begin{longlist}[(a)]
\item[(a)] Either \textup{(i)} suppose that $f$ has bounded increase at $v$, and let
$\theta=
\pr(\po(\lambda(n - 1)\sigma) > b)$ or \textup{(ii)} suppose that
$f$ has strongly bounded increase at~$v$, and let
$\theta= (n - 1) \pr(\po(\lambda d \sigma) > b)$.
Then
%
%
\begin{equation}
\label{eqnstretch1} \pr\bigl( E \cap\bigl\{f(X_t) \leq a\mbox{ for some
} t \in[t_1,t_1 + \tau] \bigr\} \bigr) \leq\biggl
\lceil\frac{\tau}{\sigma}\biggr\rceil\bigl(g(a + b)+ \theta\bigr)
\end{equation}
and
%
%
\begin{equation}
\label{eqnstretch2} \pr\bigl( E \cap\bigl\{ f(X_t) \geq a + b\mbox{ for
some } t \in[t_1,t_1 + \tau] \bigr\} \bigr) \leq\biggl
\lceil\frac{\tau}{\sigma}\biggr\rceil\bigl( h(a)+ \theta\bigr).
\end{equation}
\item[(b)] Suppose that $f$ has bounded increase via $v$, and let
$\theta= \pr(\po(\lambda d (n - 1)\sigma/2) > b)$.
Then
%
%
\begin{equation}
\label{eqnstretch3} \quad\pr\bigl( E \cap\bigl\{ f(X_t)\leq a\mbox{ for some
} t \in[t_1,t_1 + \tau] \bigr\} \bigr) \leq\biggl
\lceil\frac{\tau}{\sigma}\biggr\rceil\bigl(g(a + 2b)+ \theta\bigr)
\end{equation}
and
%
%
\begin{equation}
\label{eqnstretch4} \quad\pr\bigl( E \cap\bigl\{f(X_t) \geq a + 2b\mbox{ for
some } t \in[t_1,t_1 + \tau] \bigr\} \bigr) \leq\biggl
\lceil\frac{\tau}{\sigma}\biggr\rceil\bigl(h(a)+ \theta\bigr).
\end{equation}
\end{longlist}
\end{lemma}

\begin{pf}
We may assume that $\frac{\tau}{\sigma}$ is a positive integer $j$
by considering replacing~$\sigma$ by
$\sigma'= \tau/ \lceil\frac{\tau}{\sigma} \rceil$ (note that
$0<\sigma' \leq\sigma$ and
$\frac{\tau}{\sigma'}= \lceil\frac{\tau}{\sigma} \rceil$).

Consider first the case (a)(i), when $f$ has bounded increase at $v$.
Note that the union of the $j$
intervals $I_r = [t_1 + (r-1)\sigma, t_1+ r\sigma]$ for $r=1,\ldots,j$
is $[t_1,t_1+\tau]$.
Let $A_r$ denote the event that there are ${>}b$ arrivals for node $v$
in the interval $I_r$,
so that $\pr(A_r)= \pr(\po(\lambda(n-1) \sigma) > b)=\theta$.
Observe that, given $t \in I_r$, if $f(X_t) \leq a$ and $A_r$ fails
(so there are at most $b$ arrivals during $I_r$), then $f(X_{t_1 + r
\sigma}) \leq a+b$.
Thus
\begin{eqnarray*}
&& E \cap\bigl\{f(X_t) \leq a \mbox{ for some } t \in[t_1,t_1+
\tau]\bigr\}
\\
&&\qquad \subseteq E \cap\Biggl\{ \Biggl(\bigcup_{r=1}^{j}
\bigl\{f(X_{t_1 + r\sigma}) \leq a+b\bigr\} \Biggr) \cup\Biggl( \bigcup
_{r=1}^{j} A_r \Biggr) \Biggr\},
\end{eqnarray*}
and~(\ref{eqnstretch1}) follows.
Similarly
\begin{eqnarray*}
&& E \cap\bigl\{f(X_t) \geq a+b \mbox{ for some } t
\in[t_1,t_1+\tau]\bigr\}
\\
&&\qquad \subseteq E \cap\Biggl\{ \Biggl(\bigcup_{r=1}^{j}
\bigl\{f(X_{t_1 + (r-1)\sigma}) \geq a \bigr\} \Biggr) \cup\Biggl
(\bigcup
_{r=1}^{j} A_r \Biggr) \Biggr\},
\end{eqnarray*}
and~(\ref{eqnstretch2}) follows.

To handle the case (a)(ii) when $f$ has strongly bounded increase at $v$,
note that the arrival process onto any given link $vu$ is stochastically
dominated by a Poisson process with rate
\[
(n-2) \lambda{ {(n-2)^d -(n-3)^d} \over{(n-2)^d}} \le\lambda d.
\]
[Here we used the inequality $(1-x)^d \geq1-dx$ for $0 \leq x \leq1$.]
Thus if $B_r$ denotes the event that there are ${>}b$ arrivals in the
interval $I_r$
that are routed on some link $vu$, $u \neq v$, then
\[
\pr(B_r) \leq(n - 1) \pr\bigl(\po(\lambda d \sigma) > b\bigr);
\]
and we can complete the proof as above, replacing events $A_r$ with
events $B_r$.

Finally, in the case (b) the arrival process onto links with $v$ as
the intermediate node is stochastically dominated by a superposition of
${{n-1} \choose2}$ independent Poisson processes, each with rate
\[
\lambda{ { (n-2)^d -(n-3)^d}\over{(n-2)^d}} \le{{\lambda d} \over{n-2}}.
\]
If $C_r$ denotes the event that there are ${>}b$ arrivals in the
interval $I_r$
that are routed via $v$, then
$\pr(C_r)\le\pr(\po(\lambda d(n - 1) \sigma/2) > b)$.
The rest of the proof is as above.
\end{pf}


We present one more lemma in this subsection.
Consider a continuous-time jump Markov chain $M=(M_t)_{t \ge0}$ with
countable state space $S$
and with $q$-matrix $q=(q(x,y)\dvtx  x,y \in S)$.
Under certain conditions we can compare features of its behaviour with
that of independent birth-and-death processes.
We shall need the following lemma to handle the lower bound part of
Theorems~\ref{thmfdar1}~and~\ref{thmfdar2}.

Let $N$ be a positive integer,
and let the index $j$ run over $\{1,\ldots,N\}$.
For each $j$ let $e_j$ denote the $j$th unit $N$-vector,
and let $f_j$ be a function from $S$ to $\ints^+$,
and write $f(x)$ for $(f_1(x),\ldots,f_N(x))$.
Assume that the following two conditions hold:
\begin{longlist}[(ii)]
\item[(i)] for all distinct $x$ and $y$ in $S$ such that $q(x,y)>0$,
we have
$f(y)=f(x) \pm e_j$ for some $j$;
\item[(ii)] for each $x \in S$ and each $j$
\[
\sum_{y \in S\dvtx  f_j(y)=f_j(x)-1} q(x,y) = f_j(x).
\]
\end{longlist}
Now define $\lambda_j(x)$ for each $x \in S$ and each $j$ by setting
\[
\lambda_j(x) = \sum_{y \in S\dvtx  f_j(y)=f_j(x)+1} q(x,y).
\]
%

\begin{lemma} \label{lemnew}
Let $M$ be a continuous-time jump Markov chain as above.
For each $j$ let $\lambda_j >0$ be a constant.
Let $0 \leq t_1 < t_2$. For $j=1, \ldots, N$, let $W(j)= (W_t(j))_{t
\ge0}$ be independent birth-and-death
processes, where each $W(j)$ has constant birth rate $\lambda_j$ and
death rate equal to $w$ when in state $w$,
where $W_0(j)=0$ for each $j$.
Let $W=(W(j)\dvtx  j = 1, \ldots, N)$.
Let $F \subseteq S$ be such that for each $x \in F$ and each $j$
we have $\lambda_j(x) \geq\lambda_j$, and let $A$ be the event that
$M_t \in F$ for each $t \in[t_1,t_2]$.
Then for each downset $B$ in $\{0,1,\ldots\}^N$,
%
%
\begin{equation}
\label{eqnnew1} \pr\bigl(\bigl\{f(M_{t_2}) \in B\bigr\} \cap A\bigr)
\leq\pr(W_{t_2 -t_1} \in B).
\end{equation}
Now let $n_j$ be a given positive integer for each $j=1, \ldots, N$.
Let $\tilde{W} = (\tilde{W}_t)_{t \ge0}$, where $\tilde
{W}_t=(\tilde{W_t}(j)\dvtx  j=1, \ldots, N)$ and each $\tilde{W}(j)=
(\tilde{W}_t(j))_{t \ge0}$ is like $W(j)$
except that $\tilde{W}(j)$ has upper population limit $n_j$.
Let $\tilde{F} \subseteq S$ be such that, for each $x \in\tilde{F}$
and each $j=1,\ldots,N$,
if $f_j(x) < n_j$, then $\lambda_j(x) \geq\lambda_j$.
Let $\tilde{A}$ be the event that $M_t \in\tilde{F}$ for each $t \in
[t_1,t_2]$.
Then for each downset $B$ in $\{0,1,\ldots\}^N$,
%
%
\begin{equation}
\label{eqnnew2} \pr\bigl(\bigl\{f(M_{t_2}) \in B\bigr\} \cap\tilde{A}
\bigr) \leq\pr(\tilde{W}_{t_2
-t_1} \in B).
\end{equation}
\end{lemma}

\begin{pf}
Let us prove~(\ref{eqnnew1}), the first part of the lemma:
the second part, with population limits, may be proved similarly.

Let $W'_t(j) = W_{t-t_1}(j)$, and similarly for $W'_t$ and $W'$.
Let $m_0 \in F$, and condition on $M_{t_1}=m_0$.
Then we may assume that $\lambda_j(x) \geq\lambda_j$ for
each $x \in S$, since the values $\lambda_j(x)$ for $x \notin F$
are irrelevant, and then we may ignore the event $A$. But now we can
couple $M$ with $W'$ in
such a way that, for each $j$, every arrival in $W'(j)$ is matched by
an increment in $f_j$.
Also, for each $j$, whenever $f_j (M_t) = W'_t(j)$,
every event decreasing $f_j$ can be matched by a departure in $W'(j)$.
Since $f_j(M_{t_1}) \ge0$ for all $j$, under the coupling, $f_j (M_t)
\ge W'_t(j)$ for
each $j=1, \ldots, N$ and $t \in[t_1,t_2]$, and it follows in
particular that
\[
\pr\bigl(f(M_{t_2}) \in B \mid M_{t_1}=m_0\bigr)
\leq\pr\bigl(W'_{t_2} \in B\bigr)= \pr(W_{t_2 -t_1}
\in B).
\]
Inequality~(\ref{eqnnew1}) now follows
since this is true for each $m_0 \in F$.
\end{pf}


\subsection{PASTA}\label{subpasta}

We shall need information on the behaviour of our routing systems at
arrival times of calls, and sometimes we will need to use the following
nonasymptotic version of the
PASTA principle (``Poisson arrivals see time averages'').

Let $M=(M_t)_{t \geq0}$ be a Markov process with state space $S$, let
$(N_t)_{t \geq0}$
be a~Poisson ``arrival'' process with constant rate $\lambda$, and
assume that for each
$s >0$, $M_s$ and the process $(N_t - N_s)_{t \geq s}$ are independent.
Thus we have the ``lack of anticipation property'' that for each time $s$
the future arrivals are independent of the process up to time $s$.
Let $f$ be a bounded real-valued function on $S$.

Let $0 \leq a <b$ be fixed.
Let $V$ be the sum of the values $f(M_{t-})$ over the arrival times $t$
in $[a,b]$.
We are interested in $\mathbb EV$.

%
\begin{lemma} \label{lempasta1}
Let $\alpha= \inf_{t \in[a,b]} \mathbb E[f(M_t)]$ and
$\beta= \sup_{t \in[a,b]} \mathbb E[f(M_t)]$. Then
\[
\alpha\lambda(b-a) \leq\mathbb EV \leq\beta\lambda(b-a),
\]
and in particular, if $M_t$ is stationary, then $ \mathbb EV = \lambda
(b-a)\cdot\mathbb E[f(M_a)]$.
\end{lemma}

For example, consider a simple queuing system in equilibrium,
as in the work of Anagnostopoulos et al.~\cite{aku} on routing random calls.
Here we have an $M/M/B/B$ queue, where the Poisson arrivals have rate
$\lambda$,
service times are exponential, and there are $B$ servers;
and further where there can be at most $B$ customers in the system.
Let $P_1$ be the probability that there are $B$ customers in the system.
Then the expected number of customers lost in a unit time interval is
$\lambda P_1$.

To deduce this from the lemma above,
take $M_t$ as the number of customers at time $t$, and let $f(B)=1$ and
$f(x)=0$ if $x \neq B$,
so that $V$ is the number of customers lost.
(It is not true that the probability that a customer is lost is
$(1-e^{-\lambda}) P_1$,
as stated in the proof of Theorem~6 of~\cite{aku}.)

\begin{pf*}{Proof of Lemma~\ref{lempasta1}}
Let $A$ be the number of arrivals in $[a,b]$.
For each $k=1,2,\ldots$ on the event that $A \geq k$, let $T_k$ be the
$k$th last arrival time
in $[a,b]$ and let $V_k = f(M_{T_k-})$;
and otherwise let $T_k=-1$ say and let $V_k=0$. Then \mbox{$V = \sum_{k \geq
1} V_k$}.

First let us consider $k=1$: we shall show that
%
%
\begin{equation}
\label{eqnpasta1} \mathbb EV_1= \int_{0}^{b-a}
\mathbb E\bigl[f(M_{t})\bigr] \lambda e^{-\lambda
(b-t)} \,dt.
\end{equation}
To prove this result, let $c = b-a$ and
for each $n=1,2,\ldots$ let
\[
I_n= \biggl\lceil\frac{n (T_1-a)}{c} \biggr\rceil-1,
\]
and note that $T_1-\frac{c}{n} \leq a+\frac{cI_n}{n} < T_1$.
Let
\[
Y_n = f(M_{a+(cI_n/n)}) = \sum_{i=0}^{n-1}
\ind_{\{a+ (ci/n) <T_1 \leq a+ (c(i+1)/n) \}} f(M_{a+(ci/n)}).
\]
Then $Y_n \to V_1$ a.s.,
and so $\mathbb EY_n \to\mathbb EV_1$ by dominated convergence.
Also, crucially, the random variables
$\ind_{\{a+(ci/n) <T_1 \leq a+ (c(i+1)/n)\}}$ and
$f(M_{a+(ci/n)})$ are\vspace*{1pt} independent for each~$i$.
Hence, since $b-T_1$ has probability density $\lambda e^{-\lambda t}$
for $0\le t\le c$,
\begin{eqnarray*}
\mathbb EY_n & = & \sum_{i=0}^{n-1}
\pr\biggl(a+ \frac{ci}{n} <T_1 \leq a+ \frac
{c(i+1)}{n}
\biggr) \mathbb E\bigl[f(M_{a+(ci/n)})\bigr]
\\
& = & \sum_{i=0}^{n-1} \pr\biggl( c
\biggl(1- \frac{i+1}{n}\biggr) \leq b-T_1 < c\biggl(1-
\frac{i}{n}\biggr) \biggr) \mathbb E\bigl[f(M_{a+(ci/n)})\bigr]
\\
& = & \bigl(e^{\lambda c/n} -1\bigr) \sum_{i=0}^{n-1}
e^{-\lambda c(1 - (i/n))} \mathbb E\bigl[f(M_{a+(ci/n)})\bigr]
\\
& \sim& \lambda\frac{c}{n} \sum_{i=0}^{n-1}
\mathbb E\bigl[f(M_{a+(ci/n)})\bigr] e^{-\lambda(b- (a+(ci/n)))}
\\
& \to& \int_{0}^{b-a} \mathbb E
\bigl[f(M_{t})\bigr] \lambda e^{-\lambda(b-t)} \,dt\qquad\mbox{as } n \to
\infty
\end{eqnarray*}
since
$\mathbb E[f(M_{t})]$ is continuous as a function of $t$.
This establishes~(\ref{eqnpasta1}).

Now consider general $k \geq1$. Denote the probability density function
of $T_k$ on $[a,b]$ by $g_k(t)$.
Then just as for~(\ref{eqnpasta1}) we have
%
%
\begin{equation}
\label{eqnpasta2} \mathbb EV_k= \int_{0}^{b-a}
\mathbb E\bigl[f(M_{t})\bigr] g_k(b-t) \,dt.
\end{equation}
But $V = \sum_{k \geq1} V_k$ and
\[
\sum_{k \geq1} \int_{0}^{b-a}
g_k(b-t) \,dt = \sum_{k \geq1} \pr(A \geq k) =
\mathbb EA = \lambda(b-a),
\]
and the lemma follows.
\end{pf*}


\section{Saturated links and failure probability}\label{sectotal}

In this section we give lemmas specific to the network setting.
We give upper and lower bounds on the conditional blocking probability
of a call,
upper and lower bounds on the total number of active calls for a node
$v$, and upper bounds
on the number of saturated links incident on~$v$ over long periods of time.
All the results are valid for any GDAR algorithm.

Recall that $X_t = (X_t (e,w)\dvtx  e \in E, w \in V \setminus e)$ denotes
the load vector at time~$t$.
For each time $t$ we let $\mathcal F_t$ denote the $\sigma$-field
generated by $(X_s\dvtx  s \le t)$
(i.e., the $\sigma$-field of events up to and including time $t$).
Given a stopping time $T$ with respect to this filtration,
we let $\mathcal F_T$ denote the $\sigma$-field of all events up to
and including time~$T$,
and let $\mathcal F_{T-}$ denote the $\sigma$-field of events strictly
before~$T$.

First we consider the failure probability of a call.
Recall that
for $k=1,2,\ldots$ the call $Z_k$ arrives at time $T_k$.
For each $k$
and each node $v$, for brevity let $d_{(k)}(v)$ denote $S_{T_k-}(\mbox
{at } v) $, the number of full links
at $v$ when the call $Z_k$ arrives.
The next lemma is central to our results.

%
\begin{lemma} \label{lemblocking}
For each $k=1,2,\ldots$
%
%
\begin{equation}
\label{eqblock1} \pr(Z_k \mbox{ fails}\mid Z_k,\mathcal
F_{T_k-}) \leq{ \biggl(\frac{2 \max_v d_{(k)}(v)}{n-2} \biggr)}^d
\end{equation}
and
%
%
\begin{equation}
\label{eqblock2} \pr(Z_k \mbox{ fails}\mid\mathcal
F_{T_k-}) \le\frac{2^{d+1}}{n} \sum_{v \in V}
{ \biggl(\frac{d_{(k)}(v)}{n-2} \biggr)}^d;
\end{equation}
also, assuming that $n \geq4$,
%
%
\begin{equation}
\label{eqblock3} \pr(Z_k \mbox{ fails} \mid\mathcal
F_{T_k-}) \ge\frac{1}{2n} \sum_{v \in V}
{ \biggl(\frac{d_{(k)}(v)}{n-2} \biggr)}^d.
\end{equation}
\end{lemma}

\begin{pf}
Conditional on the event that, when a call arrives, the sum of the
numbers of saturated links at the ends
of the call is $s$, the probability it fails is at most $(\frac{s}{n-2})^d$.
Thus
\begin{eqnarray*}
\pr\bigl(Z_k \mbox{ fails}\mid Z_k=\{u,v\}, \mathcal
F_{T_k-}\bigr) & \le& { \biggl( \frac{d_{(k)}(u) + d_{(k)}(v)}{n-2}
\biggr)
}^d
\\
& \le& \biggl(\frac{2 \max_w d_{(k)}(w)}{n-2} \biggr)^d,
\end{eqnarray*}
which gives~(\ref{eqblock1}). Similarly,
\begin{eqnarray*}
\pr(Z_k \mbox{ fails}\mid\mathcal F_{T_k-}) & \le&
\frac{1}{{n \choose
2}} \sum_{u \neq v} { \biggl(
\frac{d_{(k)}(u) + d_{(k)}(v)}{n-2} \biggr) }^d
\\
& \le& \frac{2^{d-1}}{{n \choose2}} \sum_{u \neq v} \biggl(
\frac{d_{(k)}(u)^d +d_{(k)}(v)^d}{n-2} \biggr)^d,
\end{eqnarray*}
and~(\ref{eqblock2}) follows.\vspace*{1pt}
[For the second inequality we used the fact that $f(x)=x^d$ is convex
for $x>0$, and so
$(x+y)^d \leq2^{d-1}(x^d +y^d)$ for $x,y>0$.]

On the other hand,
\[
\pr(Z_k \mbox{ fails}\mid\mathcal F_{T_k-}) \ge
\frac{1}{2{n \choose2}} \sum_{v \in V} \sum
_{u \neq v} \biggl( \frac{d_{(k)} (v)- \ind_{X_{T_k-}(vu)=D}}{n-2}
\biggr)^d.
\]
But for each $v \in V$,
\begin{eqnarray*}
&& \sum_{u \neq v} \bigl(d_{(k)}(v)-
\ind_{X_{T_k-}(vu)=D}\bigr)^d
\\
&&\qquad = \bigl(n-1-d_{(k)}(v)\bigr) d_{(k)}(v)^d +
d_{(k)}(v) \bigl(d_{(k)}(v)-1\bigr)^d
\\
&&\qquad  \geq \frac{1}2 (n-1) d_{(k)}(v)^d
\end{eqnarray*}
for $n \geq4$.
[To see this, consider separately the case $d_{(k)}(v) \leq\frac{n-1}{2}$,
when the first term suffices; and the case $d_{(k)}(v)=x \geq\frac{n}{2}$,
when $(x-1)^d \geq x^d (1-d/x) \geq\frac{1}2 x^d$.]
Hence
\begin{eqnarray*}
\pr(Z_k \mbox{ fails}\mid\mathcal F_{T_k-}) & \ge&
\frac{1}{2{n \choose2}} \sum_{v \in V}\frac{ 1/2 (n-1)
d_{(k)}(v)^d}{(n-2)^d}
\\
&=& \frac{1}{2n} \sum_{v \in V} { \biggl(
\frac{d_{(k)}(v)}{n-2} \biggr)}^d,
\end{eqnarray*}
and~(\ref{eqblock3}) follows.
\end{pf}

To obtain our estimates for the total number of active calls for a node $v$,
and upper bounds on the number of saturated links incident on $v$,
we compare the process~$X$ to a much simpler dominating process $\tilde
{X} = (\tilde{X}_t)_{t \ge0}$
which also has state space ${\mathcal X}$ and satisfies $\tilde
{X}_0=X_0$ and evolves as follows.
The
edges $e=\{u, v\}$ in $E(K_n)$ receive independent rate
$\lambda$ Poisson arrival streams of calls; each link $uv$ has
infinite capacity, and each call throughout its duration occupies
$d$ two-link routes chosen uniformly at random with replacement.
(If a route is chosen more than once by a given call, then the call will
still be counted only once on the corresponding two links.) All call
durations are unit mean exponentials independent of one another and
of the arrivals and choices processes, and whenever a call is
completed, it frees all the links it has been occupying.

As for process $X$, for each edge $e$ in $E(K_n)$ and each node $w
\notin e$, we let $\tilde{X}_t(e,w)$ denote
the number of calls between the end nodes of $e$ routed via $w$ in
progress at time $t$;
also, let $\tilde{X}_t(vu)$ denote the load of link $vu$,
let $\tilde{X}_t (e)$ denote the number of calls in progress between
the end nodes of $e$ at time $t$,
and let $\tilde{X}_t(v)$ denote the number of calls with one end $v$
in progress at time $t$.
[Note that, in contrast to process $X$, here it is not necessarily
the case that
$\tilde{X}_t(e)$ equals $\sum_{w \notin e} \tilde{X}_t (e,w)$,
and it is not the case that $\tilde{X}_t (v)$ equals $\sum_{u \neq v}
\tilde{X}_t(vu)$;
this is because a single call is allowed to occupy more than one route.]
Note further that the process
$(\tilde{X}_t(e)\dvtx  e \in E(K_n))$ is itself
Markov, since the capacities are infinite, and so no calls get rejected.
It has a unique equilibrium distribution,
and in equilibrium the $\tilde{X}_t(e)$ are all independent
$\po(\lambda)$ random variables. Thus in equilibrium, the total
number $\llVert\tilde{X}_t \rrVert_1$ of ongoing calls at
time $t$ is
$\po(\lambda{n \choose2})$; and, for each~$v$, the total number
$\tilde{X}_t (v)$ of ongoing calls with one end $v$ is
$\po(\lambda(n-1))$.

We shall use $T(v)$ to denote the time that the last of the $X_0(v)$
initial calls with one end $v$ departs. Also, we let $T = \max_{v \in
V} T(v)$,
the time when the last of the initial $\llVert X_0\rrVert_1$ calls depart.
As was mentioned earlier, if initially there are many calls, then the
system needs a
``burn-in'' period to reduce ``congestion'' measures such as the number
of full links. The system will have
``lost'' the memory of the bad initial state once all the initial calls
are completed. For this reason,
for various events $A$ we shall give an upper bound on
$\pr(A \cap\{ T \leq t\})$. We may later obtain an upper bound on
$\pr(A)$ using
%
%
\begin{equation}
\label{eqnTt1} \pr(A) \leq\pr\bigl(A \cap\{ T \leq t\}\bigr) + \pr(T>t),
\end{equation}
and noting that
%
%
\begin{equation}
\label{eqnTt2} \pr(T > t) \le\mathbb E\llVert X_0 \rrVert
_1 e^{-t}.
\end{equation}
To see why~(\ref{eqnTt2}) holds, temporarily let $I_t$ denote the
number of initial calls
surviving to time $t$, and observe that
\[
\pr(T > t) = \pr(I_t>0) \leq\mathbb EI_t = \mathbb E
\llVert X_0 \rrVert_1 e^{-t}.
\]
We shall always be interested in link capacities $D(n)$ which grow
slowly with $n$,
and which in particular satisfy $D(n)=o(n)$. Thus always $\llVert
X_0 \rrVert_1
=o (n^3)$,
and so~(\ref{eqnTt2}) gives
%
%
\begin{equation}
\label{eqnTt3a} \pr(T > t) = o\bigl(n^3\bigr) \cdot e^{-t}.
\end{equation}
If $X_0$ is stochastically at most the equilibrium distribution for
$\tilde{X}$,
we let $\tilde{T}=0$ a.s. and otherwise let $\tilde{T}=T$.

The next lemma shows that for any node $v \in V$,
the number $X_t(v)$ of calls at $v$ is unlikely to deviate far above
$\lambda(n-1)$
once the initial calls have gone.

%
\begin{lemma} \label{lemcalls1}
Let $0<\delta< 1$, let $n$ be a positive integer, and
let $A_t$ be the event that
$X_t(v)\ge(1+\delta)\lambda(n-1)$ for some vertex $v$.
Then for all times $t_1 \geq0$ and $t_2 \ge t_1$,
%
%
\begin{equation}
\label{calls1b} \pr\bigl(A_{t_2} \cap\{\tilde{T} \leq t_1
\}\bigr)\leq n e^{-(1/3)\delta^2 \lambda(n-1)}.
\end{equation}
\end{lemma}

Note that the value of $D$ is not relevant here.

\begin{pf*}{Proof of Lemma~\ref{lemcalls1}}
Let $\tilde{Y} = (\tilde{Y}_t)_{t \ge0}$, with
$\tilde{Y}_t= (\tilde{Y}_t(e,w)\dvtx  e \in E(K_n),\break  w \notin e)$,
be a Markov process with the same $q$-matrix as $(\tilde{X}_t)$ but in
equilibrium.
Observe\vspace*{2pt} that the equilibrium distribution for $(X_t)$ is stochastically
at most the
distribution for $\tilde{Y}_t$.
We couple $(X_t)$, $(\tilde{X}_t)$, and
$(\tilde{Y}_t)$ as follows.
We assume that $X_0 = \tilde{X}_0$, and further if $X_0$ is
stochastically at most the equilibrium distribution~$\tilde{Y}_0$, then $X_0 = \tilde{X}_0 \leq\tilde{Y}_0$.
All subsequent arrival and potential departure times of new calls are the
same for the three processes, except that the departures of calls that
were not
accepted due to none of their chosen routes being available in $(X_t)$
are ignored in that process.
Additionally,\vspace*{1.5pt} every one of the $\llVert X_0 \rrVert_1$
initial calls
in $(X_t)$ is coupled with a corresponding initial call in
$(\tilde{X}_t)$ and in $(\tilde{Y}_t)$, and the paired calls have the same
departure times.

Since\vspace*{2pt} all calls are accepted in $\tilde{X}$ and in $\tilde{Y}$,
under the coupling, for each node $v$ and time $t$,
on the event $\tilde{T} \leq t$ we have
%
%
\begin{equation}
\label{eqncouple2} X_t(v)\le\tilde{X}_t(v)\le
\tilde{Y}_t (v).
\end{equation}
But $\tilde{Y}_t (v)$ is a Poisson random variable
with mean $\lambda(n-1)$, and so by the concentration
inequality~(\ref{binpo1}), we have, for each $v$, and all $t_2 \ge t_1$,
\[
\pr\bigl( \bigl\{X_{t_2}(v)\ge(1 + \delta)\lambda(n-1)\bigr\} \cap\{
\tilde{T} \leq t_1\} \bigr) \leq e^{-(1/3)\delta^2 \lambda(n-1)}.
\]
Now~(\ref{calls1b}) follows by summing the above bound over all $v$.
\end{pf*}

We will now use the above result to show that, after a burn-in period,
we are unlikely to observe large deviations of $X_t(v)$ above
$\lambda(n-1)$ for any node $v$ even during very long time intervals.
Recall the notation $p_D(\mu)$ introduced in~(\ref{eqnpobound2}).
%

\begin{lemma} \label{lemcalls3}
Given $0<\delta< 1$, there exists a constant $\beta=\beta(\delta)
>0$ such that
the following holds.
Let the capacity $D=D(n)=o(n)$.
Let $\kappa> 0$, and let $\tilde{t}_0 = (\kappa+ 3) \ln n$.
If $X_0$ is stochastically at most the equilibrium distribution, let
$t_0=0$, and otherwise let
$t_0=\tilde{t}_0$.
Let $C_t$ denote the event that $X_t (v) > (1 + \delta) \lambda(n-1)$
for some vertex~$v$.
Then as $n \to\infty$, for each time $t_1 \geq t_0$
\[
\pr\bigl(C_t \mbox{ holds for some } t \in\bigl[t_1,t_1+e^{\beta
n}
\bigr] \bigr) = o\bigl(n^{-\kappa}\bigr).
\]
\end{lemma}

\begin{pf}
Let $C'_t$ denote the event that $X_t (v) >(1+ \delta/2) \lambda
(n-1)$ for some vertex~$v$.
By Lemma~\ref{lemcalls1},
there exists a constant $\gamma>0$ such that for each time $t \geq t_0$,
\[
\pr\bigl(C'_t \cap\{\tilde{T} \leq t_0\}
\bigr) \le2e^{-\gamma n}.
\]
We may assume that $\gamma\leq\delta/12$.
Let $\beta=\gamma/3$. Let $v \in V$, and let $f(X_t)=X_t(v)$, which
has bounded increase at $v$.
We now apply inequality~(\ref{eqnstretch2}) in Lemma~\ref
{lemstretch}, with
$a = (1+\delta/2) \lambda(n-1)$, $b = (\delta/2) \lambda(n-1)$,
$\tau= e^{\beta n}$, $\sigma= \delta/4$, and $E$ the event $\{\tilde
{T} \leq t_0\}$.
Also, let $\theta= \pr(\po (\lambda(n-1)\delta/4) > \lambda
(n-1) \delta/2)$.
Thus for all positive integers $n$ and all times $t_1 \geq t_0$,
we have
\begin{eqnarray*}
&& \pr\bigl( \bigl\{ X_t(v) >(1 + \delta) \lambda(n - 1) \mbox{ for
some } t \in\bigl[t_1,t_1 + e^{\beta n}\bigr]
\bigr\} \cap\{\tilde{T} \leq t_0 \} \bigr)
\\
&&\qquad \leq\bigl((4/ \delta) e^{\beta n} +1 \bigr) \bigl(2e^{-\gamma n} +
\theta\bigr).
\end{eqnarray*}
Also, by~(\ref{binpo1}), $\theta\leq e^{-(n-1)\delta/12} =
O(e^{-\gamma n})$. Hence,
summing over the $n$ nodes in $V$, we obtain
\[
\pr\bigl( \bigl\{ C_t \mbox{ for some } t \in\bigl[t_1,t_1+e^{\beta n}
\bigr] \bigr\} \cap\{\tilde{T} \leq t_0 \} \bigr) = o
\bigl(e^{-\beta n}\bigr).
\]
We may now use~(\ref{eqnTt1}) and~(\ref{eqnTt2}) to complete the proof,
noting that always $\llVert X_0 \rrVert_1 = O(n^2D) =o(n^3)$.
\end{pf}

To end this section we
shall upper bound the number of saturated links around any given node
in the following lemma.
Observe from~(\ref{eqnpobound1}) that if we have $\delta> 0$ and
$D=D(n) \to\infty$,
then for $n$ sufficiently large we may, for example, take $k$ as $\delta
n$ in the lemma.
%

\begin{lemma} \label{lemsat}
Let $n$ and $D$ be positive integers, and let $k \geq4 p_D(d\lambda) (n-1)$.
Then for each $t \geq0$,
%
%
\begin{equation}
\label{eqnlinks1} \pr\bigl( \bigl\{S_{t}(\mbox{at } v) \ge k \bigr\}
\cap\{\tilde{T} \leq t\} \bigr) \leq2 \exp{ \biggl(-{{k} \over{16 d^2 D}}
\biggr)}
\end{equation}
and
%
%
\begin{equation}
\label{eqnlinks2} \pr\bigl( \bigl\{S_{t}^{}(\mbox{via } v)
\ge k \bigr\} \cap\{\tilde{T} \leq t\} \bigr) \leq2 \exp{ \biggl(-
{{k } \over{64 D}} \biggr)}.
\end{equation}
\end{lemma}

\begin{pf}
We use the coupling of the three processes $(X_t)$, $(\tilde{X}_t)$
and $(\tilde{Y}_t)$ described in the proof of Lemma~\ref{lemcalls1}.
Consider a link $vu$ (where $u \neq v$) and a time $t$:
under the coupling, on the event that $\tilde{T} \leq t$,
%
%
\begin{equation}
\label{eqncouple1} X_t(vu)\le\tilde{X}_t(vu)\le
\tilde{Y}_t (vu).
\end{equation}
We can thus work mostly with the stationary dominating process $(\tilde
{Y}_t)$, where
we bound expectations and use concentration inequalities.

Let $v \in V$ be a node. Note that for each $u \neq v$,
the load $\tilde{Y}_t (vu)$ of link $vu$ is a
Poisson random variable with mean
\[
\lambda(n-2){{(n-2)^d-(n-3)^d} \over{(n-2)^d}} \le d \lambda.
\]
We adapt some more notation introduced earlier for $(X_t)$ to $(\tilde
{Y}_t)$ in the natural way.
Thus we write ${\tilde{\mathcal S}}_{t}(\mbox{at } v) $ to denote the
set of links $vw$ for calls at
$v$ that have load at least $D$ at time $t$ in
$(\tilde{Y}_t)$, and we write $\tilde{S}_{t}(\mbox{at } v) =
\llvert{\tilde{\mathcal S}}_{t}(\mbox{at } v) \rrvert$.
Also, for $w \in V$, ${\tilde{\mathcal S}}_{t}(\mbox{via } w) $
denotes the set of links $uw$
for calls at some node $u$, and routed via~$w$, that have load at least
$D$ at time $t$
in $(\tilde{Y}_t)$ and $\tilde{S}_{t}(\mbox{via } w) =\llvert{\tilde
{\mathcal S}}_{t}(\mbox{via } w) \rrvert$.
Then $\mathbb E[\tilde{S}_{t}(\mbox{at } v) ] \le(n-1)p_{D}(d\lambda
)$ and
$\mathbb E[\tilde{S}_{t}(\mbox{via } w) ] \le(n-1)p_{D}(d\lambda)$
for all times \mbox{$t \ge0$}.

For a given $v \in V$, we may think of the loads $\tilde{Y}_t (vu)$ of
links $vu$
for $u \neq v$ as being determined by a family of $(n-1)(n-2)^d$
independent Poisson random variables each with mean
$\lambda/(n-2)^d$ [corresponding to $n-1$ choices of the other end
node $w$
and $(n-2)^d$ choices of $d$ routes for a call with end nodes $v$ and $w$],
and so there is strong concentration of measure.
Note that the median $m(v)$ of $\tilde{S}_{t}(\mbox{at } v) $ is at
most $2(n-1)p_{D}(d\lambda)$.
We can use Talagrand's inequality Lemma~\ref{lemtalagrand}, with
$c=d$ and
$r=D$. This gives, for all $z \geq0$,
\[
\pr\bigl( \tilde{S}_{t}(\mbox{at } v) \ge m (v)+z\bigr)\le2 \exp{
\biggl(-{{z^2} \over{4d^2 D(m (v) +z)}} \biggr)}.
\]
Now take $z \ge2(n-1)p_{D}(d\lambda) \ge m (v)$, so that
%
%
\begin{equation}
\label{eqntildeS} \pr\bigl( \tilde{S}_{t}(\mbox{at } v) \ge2z\bigr)\le
2 \exp{ \biggl(-{z \over{8d^2 D}} \biggr)}.
\end{equation}
Similarly, given $w \in V$, the loads $\tilde{Y}_t (uw)$ of links $uw$ for
$u \neq w$ may be determined by a family of ${{n-1} \choose2}
[(n-2)^d-(n-3)^d]$
independent random variables $\po(\lambda/(n-2)^d)$ (corresponding to calls
for all possible pairs of distinct nodes $v,u \in V \setminus\{w\}$
choosing a route via node $w$).
Applying Talagrand's inequality with $c =2$ and $r=D$,
we have, for $t \ge0$ and $z \ge2(n-1)p_{D}(d\lambda)$,
\[
\pr\bigl( \tilde{S}_{t}(\mbox{via } w) \ge2z\bigr)\le2 \exp{
\biggl(-{z \over{32 D}} \biggr)}.
\]
But $X_t (vu) \le\tilde{Y}_t (vu)$ on the event that $\tilde{T} \le t$
[as we noted in~(\ref{eqncouple1})],
and we deduce that inequalities~(\ref{eqnlinks1})
and~(\ref{eqnlinks2}) hold.
\end{pf}


\section{Proof of Theorems~\texorpdfstring{\protect\ref{thmfdar1}}{1.1} 
and~\texorpdfstring{\protect\ref{thmfdar2}}{1.2}}\label{secproof1}

Let us recall the rough story.
The number $S_{t}(\mbox{at } v) $ of saturated links at a node $v$ has
expected value $n^{1-\alpha+ o(1)}$,
and by Lemma~\ref{lemblocking} the probability $p$ that a call with
one end $v$ fails is roughly $\mathbb E[(S_{t}(\mbox{at } v) /n)^d]$.
There is a change of behaviour at $\alpha=1$. If $0<\alpha<1$, then
$S_{t}(\mbox{at } v) $ is concentrated and
$\mathbb E[S_{t}(\mbox{at } v) ^d] = n^{(1-\alpha)d + o(1)}$ and $p$
is $n^{-\alpha d + o(1)}$.
The expected number of arrivals in an interval of length $n^K$ is about
$n^{K+2}$,
and $n^{K+2} p = n^{K+2 - \alpha d +o(1)}$, so $\alpha$ is $K$-good
when $K+2 - \alpha d <0$,
and $\alpha$ is \mbox{$K$-}bad when $K+2 - \alpha d >0$.
When $\alpha\geq1$, then $\mathbb E[S_{t}(\mbox{at } v) ^d] \sim
\mathbb E[S_{t}(\mbox{at } v) ]$ and
$p = n^{1-\alpha-d +o(1)}$, and again we see when $\alpha$ is
$K$-good by looking at $n^{K+2} p$.

Note that the case $\alpha= 1$ is covered by our proofs: we show that
$\alpha$ is $K$-good if $K < d-2$, and $\alpha$ is $K$-bad if $K > d-2$.

\subsection{Upper bounds: Showing \texorpdfstring{$\alpha$}{$alpha$} is $K$-good}\label{upper1}

Here our aim is to prove that, for appropriate $\alpha$ and $K$,
if we use any GDAR algorithm on a network with $n$ nodes,
and the link capacity $D(n)\sim\alpha\ln n / \ln\ln n$ is high enough,
then the mean number of calls that are lost over an interval of length
$n^K$ is $o(1)$ as $n \to\infty$.
To achieve this, we need to be able to show that, throughout the time interval,
there are not too many saturated (full) links in the network.

We may need to wait for a ``burn-in'' period so that any initial
congestion can dissipate,
and in fact in this case
we wait until all the initial calls have left the system.
Recall that $T$ denotes the departure time of the last of the initial
calls. Recall also that $\tilde{T}=0$ if the distribution of $X_0$ is
stochastically at most the stationary distribution, and $\tilde{T} =
T$ otherwise.
Now $\llVert X_0 \rrVert_1 \le{n \choose2} D = o(n^2 \ln
n)$ (as $n \to\infty$).
Hence, by~(\ref{eqnTt2}), for each $t>0$,
%
%
\begin{equation}
\label{eqnTt3} \pr(\tilde{T} > t) \le\pr(T > t) \leq\mathbb E\llVert
X_0 \rrVert_1 e^{-t} = o\bigl(n^2
\ln n e^{-t}\bigr).
\end{equation}

Recall that we set $t_1 \geq t_0=5 \ln n$ if $X_0$ is not
stochastically dominated by the stationary distribution, and $t_0=0$
otherwise. Let $t_2=t_1 + K \ln n$, and let $t_3=t_1 + n^K$.
Then by~(\ref{eqnTt3}), $\pr(\tilde{T} > t_1) = o(n^{-2})$
and $\pr(\tilde{T} > t_2) = o(n^{-K-2})$.
For $0 \leq t < t'$ let $N_F(t,t')$ be the number of calls that fail in
the interval $(t,t']$.
We shall show that for $j=1, 2$ we have
$\mathbb EN_F(t_j,t_{j+1}) =o(1)$, yielding $\mathbb EN_F(t_1,t_1 +
n^K) =o(1)$
as required.

For $0 \leq t < t'$, let $N_A(t,t')$ be the number of calls that arrive
in the interval $(t,t']$.
Thus $N_A(t,t') \sim\po(\lambda{n \choose2} (t'-t))$.
Let $N_1 = \lceil2 \mathbb EN_A(t_1,t_2) \rceil\sim\lambda K n^2 \ln n$,
and let
$N_2 = \lceil2 \mathbb EN_A(t_1,t_3) \rceil\sim\lambda n^{K+2}$.
Finally here note that since $D \sim\alpha\ln n/\ln\ln n$,
from~(\ref{eqnpobound2}) we have
%
%
\begin{equation}
\label{eqnpDd} p_D(d\lambda)=n^{-\alpha+o(1)}.
\end{equation}
There are two subcases, for $\alpha\leq1$ and $\alpha>1$.

Suppose first that $(K+2)/d < \alpha\leq1$ (and so $d \geq3$).
In order to upper bound the probability that a call $Z_k$ fails,
we will use Lemmas~\ref{lemsat} and~\ref{lemstretch}
to upper bound the maximum number of saturated links at any node,
and then we can use inequality~(\ref{eqblock1}) in Lemma~\ref{lemblocking}.
By inequality~(\ref{eqnlinks1}) in Lemma~\ref{lemsat} with
$k=4(n-1)p_D(d \lambda)+ \ln^3 n$,
for each $v \in V$,
\begin{eqnarray*}
&& \pr\bigl( \bigl\{ S_{t}(\mbox{at } v) \geq4(n-1)
p_D(d\lambda) +\ln^3 n \bigr\} \cap\{ \tilde{T} \leq t
\} \bigr)
\\
&&\qquad = \exp\bigl(-\Omega\bigl(\ln^3 n /D\bigr)\bigr) = \exp\bigl(-
\Omega\bigl(\ln^2 n\bigr)\bigr).
\end{eqnarray*}
For $0 \leq t < t'$, let $A_{t,t'}$ be the event that
$S_{s}(\mbox{at } v) \le4 (n-1) p_D (d \lambda) + 2 \ln^3 n$
for each vertex $v$ and each $s \in(t,t']$.
By\vspace*{1.5pt} the above inequality and Lemma \ref{lemstretch}(a)(i)
[with $\tau= n^K$, $\sigma=1/n$,
$a= 4 (n-1) p_D (d \lambda) + \ln^3 n$ and $b=\ln^3n$],
for each $t \geq0$,
\[
\pr\bigl(\overline{A_{t,t+n^K}} \cap\{\tilde{T} \le t\}\bigr) = \exp
\bigl(-\Omega\bigl(\ln^2 n\bigr)\bigr),
\]
and it follows that
for $j=1$ and 2 we have
\[
N_j \pr(\overline{A_{t_j,t_{j+1}}}) = o(1).
\]

Let $j$ be 1 or 2. List the calls arriving
after $t_j$ as $Z'_1,Z'_2,\ldots$
arriving at times
$t_j <T'_1<T'_2< \cdots.$
Since for each $k=1,2,\ldots$
\[
\bigl\{T'_k \leq t_{j+1}\bigr\} \cap
A_{t_j,t_{j+1}} 
\subseteq\bigl\{S_{T'_k-}(\mbox{at } v) \le4
(n-1) p_D (d \lambda) + 2 \ln^3 n\ \forall v
\bigr\},
\]
by
inequality~(\ref{eqblock1}) applied to these arrivals we have
\[
\pr\bigl(\bigl\{Z'_{k} \mbox{ fails}\bigr\} \cap\bigl
\{T'_k \leq t_{j+1}\bigr\} \cap
A_{t,t'} 
\bigr) \le p_0,
\]
where by~(\ref{eqnpDd})
\[
p_0 = \biggl(\frac{8(n-1)p_D(d \lambda) +4\ln^3 n}{n-2} \biggr)^d =
n^{-\alpha d +o(1)} = o\bigl(n^{-K-2}\bigr).
\]
Note also
that, if the random variable $ Y_j \sim\po(\lambda{n \choose
2}(t_{j+1} - t_1)) $, then
$ \mathbb EY_j \leq N_j /2 $, and so $\mathbb E[Y_j \ind_{Y_j>N_j}] = o(1)$.
Hence
\begin{eqnarray*}
&& \mathbb E N_F(t_j,t_{j+1})
\\
&&\qquad  = \mathbb E
\Biggl[\sum_{k=1}^{N_A(t_j,t_{j+1})} \ind_{Z'_k\ \mathrm{fails}}
\Biggr]
\\
&&\qquad =  \mathbb E\Biggl[ \sum_{k=1}^{N_A(t_j,t_{j + 1})}
\ind_{Z'_k\ \mathrm{fails}} \ind_{N_A(t_j,t_{j + 1}) \le N_j}\Biggr]
\\
&&\quad\qquad{}+ \mathbb E\Biggl[ \sum
_{k=1}^{N_A(t_j,t_{j + 1})} \ind_{Z'_k\ \mathrm{fails}}\ind_{N_A(t_j,t_{j + 1}) > N_j}
\Biggr]
\\
&&\qquad  \leq \sum_{k = 1}^{N_j} \pr\bigl(\bigl
\{Z'_k \mbox{ fails}\bigr\} \cap\bigl
\{T'_k \le t_{j+1}\bigr\}\bigr) + \mathbb E
\bigl[N_A(t_j,t_{j+1}) \ind_{N_A(t_j,t_{j+1}) > N_j}
\bigr]
\\
&&\qquad \leq \sum_{k = 1}^{N_j} \pr\bigl(\bigl
\{Z'_k \mbox{ fails}\bigr\} \cap\bigl
\{T'_k \leq t_{j+1}\bigr\} \cap
A_{t_j,t_{j+1}}\bigr) + N_j \pr(\overline{A_{t_j,t_{j+1}}})
+o(1)
\\
&&\qquad \leq N_j p_0 + o(1) = O\bigl(n^{K+2}
p_0\bigr) + o(1) = o(1).
\end{eqnarray*}
Thus $\mathbb E N_F(t_1,t_{3})=o(1)$, as required.
This completes the proof of the subcase $(K+2)/d < \alpha\leq1$.

Now consider the other subcase, where $\alpha>1$ and $\alpha> K+3-d$.
We may assume that $K \ge d-2$, and now the condition reduces simply to
$\alpha> K+3-d$.
In this subcase we need a different and somewhat more involved proof.
[Note that $p_0= \Omega(n^{-d})$ and so $n^{K+2} p_0$ may be large.]
In order to upper bound the probability that a call fails,
we will use inequality~(\ref{eqblock2}) in Lemma~\ref{lemblocking},
and to upper bound
the expected number of saturated links at a node,
we use the stationary dominating process.

Fix $j$, and as before list the calls arriving after time $t_j$ as
$Z'_1, Z'_2, \ldots$
arriving at times $t_j < T'_1< T'_2 < \cdots.$
We will show that
%
%
\begin{equation}
\label{eqn1ashow} \quad\sum_{k=1}^{\infty} \pr\bigl(
\bigl\{Z'_k \mbox{ fails}\bigr\} \cap\{\tilde{T}\le
t_j \} \cap\bigl\{T'_k \le
t_{j+1}\bigr\} \bigr) \leq(t_{j+1}-t_j)
n^{3-d-\alpha+ o(1)}.
\end{equation}
From this result we will complete the proof quickly.

Now for the details.
Note first that, since each $T'_k >t_j$ we have
\[
\{\tilde{T} \le t_j\} \cap\bigl\{T'_k
\le t_{j+1}\bigr\} \in\mathcal F_{T'_k-}.
\]
Thus, by Lemma~\ref{lemblocking} inequality~(\ref{eqblock2}),
for each $k=1,2,\ldots$
\begin{eqnarray*}
&& \pr\bigl(\bigl\{Z'_k \mbox{ fails}\bigr\} \cap\{
\tilde{T} \le t_j\} \cap\bigl\{ T'_k \le
t_{j+1}\bigr\}\bigr)
\\
&&\qquad \leq\frac{2^{d+1}}{n(n-2)^d} \sum_{v \in V}
\mathbb E\bigl[ \bigl(S_{T'_k-}(\mbox{at } v) \bigr)^d
\ind_{\tilde{T} \le t_j} \ind_{T'_k \le t_{j +1}}\bigr].
\end{eqnarray*}
Recall from the proof of Lemma~\ref{lemsat}
in Section~\ref{sectotal} that there is a coupling
involving a stationary copy $(\tilde{Y}_t)$ of the dominating process
with the following property.
On $\{\tilde{T} \le t\}$, for each $v \in V$,
the number $S_{t}(\mbox{at } v) $ of links ending in $v$ which are
saturated at time $t$ is
stochastically at most the number $\tilde{S}_{t}(\mbox{at } v) $ of
links $vu$ for $u \neq v$
such that $\tilde{Y}_t (vu) \ge D$.
Therefore, for each $k=1,2,\ldots$
\[
\pr\bigl(\bigl\{Z'_k \mbox{ fails}\bigr\} \cap\{
\tilde{T} \le t_j \} \cap\bigl\{ T'_k \le
t_{j + 1}\bigr\}\bigr) \leq\frac{2^{d+1}}{n(n - 2)^d} \sum
_{v \in V} \mathbb E\bigl[ \tilde{S}_{T'_k-}(\mbox{at }
v) ^d \ind_{T'_k \le t_{j + 1}}\bigr].
\]
Hence
\begin{eqnarray*}
&& \sum_{k=1}^{\infty} \pr\bigl(\bigl
\{Z'_k \mbox{ fails}\bigr\} \cap\{T \le
t_j \} \cap\bigl\{ T'_k \le
t_{j+1}\bigr\}\bigr)
\\
&&\qquad \leq \frac{2^{d+1}}{n(n-2)^d} \mathbb E\Biggl[\sum_{k=1}^{\infty}
\sum_{v \in
V} \tilde{S}_{T'_k-}(\mbox{at } v)
^d \ind_{T'_k \le t_{j+1}}\Biggr]
\\
&&\qquad =  \frac{2^{d+1}}{n(n-2)^d} \lambda\pmatrix{n \cr2} (t_{j+1}-t_j)
\sum_{v \in V} \mathbb E\bigl[\tilde{S}_{0}(
\mbox{at } v) ^d\bigr],
\end{eqnarray*}
where the last equality follows from the PASTA property of Lemma~\ref{lempasta1}.\vspace*{1.5pt}

Let us write $\tilde{d}(v)$ for $\tilde{S}_{0}(\mbox{at } v) $ for brevity.
We now claim that, for each $v \in V$,
%
%
\begin{equation}
\label{claimdvd} \mathbb E\bigl[\tilde{d}(v)^d\bigr] = \mathbb E\bigl[
\tilde{d}(v)\bigr] \bigl(1+o(1)\bigr)=n^{1-\alpha+o(1)}.
\end{equation}
Inequality~(\ref{eqn1ashow}) will follow immediately from
the last result and claim~(\ref{claimdvd}).

To prove the claim, consider the dominating process at time 0.
Consider first a fixed link $vw$.
The probability that a call $\{u,v\}$ uses this link is
$1-(1- \frac{1}{n-2})^d$. Thus from our earlier discussion
the load on the link has Poisson distribution with mean
$\lambda(n-1) (1-(1- \frac{1}{n-2})^d) = \lambda d + O(n^{-1})$.
It follows as in~(\ref{eqnpDd}) that
the probability that $vw$ has load at least $D$ is $n^{-\alpha+o(1)}$,
and so
\[
\mathbb E\bigl[\tilde{d}(v)\bigr] = n^{1-\alpha+o(1)}.
\]
This gives one part (the easy part) of claim~(\ref{claimdvd}).

Now fix $v \in V$, and let $u_1,\ldots,u_d$ be
distinct nodes in $V \setminus\{v\}$.
Let $N(u_i)$ be the number of live calls with one end $v$ that
have selected the link $vu_i$ but none of the links $vu_j$ for
$j \neq i$. Let $\tilde{N}$ be the number of live calls that
have selected at least two of the links $vu_i$.
Then the $N(u_i)$ are i.i.d., each is Poisson with mean at most
$\lambda d$, and
$\tilde{N}$ is Poisson with mean $O(1/n)$.

Let $x=d+\alpha$, and let $A$ be the event that $\tilde{N} \leq x$.
Note that $\pr(\bar{A}) = O(n^{-x})$ by~(\ref{eqnpobound1}) and
\[
\mathbb E \Biggl[\prod_{i=1}^{k}
\ind_{\tilde{Y}_0(vu_{i}) \ge D } \ind_{A} \Biggr] \leq\mathbb E \Biggl
[\prod
_{i=1}^{k} \ind_{N(u_i) \ge D- x} \Biggr] =
\pr\bigl(N(u_{1}) \ge D- x\bigr)^k.
\]
Also, by~(\ref{eqnpobound1}),
$\pr(N(u_{1}) \ge D-x) \leq
n^{-\alpha+ o(1)}$.
Now let $a_k$ be the number of partitions of $1,\ldots,d$ into
exactly $k$ nonempty blocks.
In the sums below the $w_{j}$ run over $V \setminus\{v\}$.
We find
\begin{eqnarray*}
\mathbb E\bigl[\tilde{d}(v)^d \ind_{A}\bigr] &=& \mathbb
E \Biggl[\prod_{j=1}^{d} \sum
_{w_{j}} \ind_{ \tilde{Y}_0(vw_{j}) \ge D } \ind_{A} \Biggr]
\\
&=& \sum_{w_{1},\ldots,w_{d}} \mathbb E\Biggl[ \prod
_{j=1}^{d} \ind_{ \tilde{Y}_0(vw_{j}) \ge D } \ind_{A}
\Biggr]
\\
&=& \sum_{k=1}^{d} a_k
(n-1)_k \mathbb E\Biggl[ \prod_{i=1}^{k}
\ind_{ \tilde{Y}_0(vu_{i}) \ge D } \ind_{A}\Biggr]
\\
& \leq& \mathbb E\bigl[\tilde{d}(v)\ind_{A}\bigr] + \sum
_{k=2}^{d} a_k n^k \pr
\bigl(N(u_{1}) \ge D-x\bigr)^k
\\
& \leq& \mathbb E\bigl[\tilde{d}(v)\bigr] + O\Biggl( \sum
_{k=2}^{d} \bigl(n^{1-\alpha+o(1)}
\bigr)^k\Biggr)
\\
& = & \mathbb E\bigl[\tilde{d}(v)\bigr] + O\bigl(n^{-2(\alpha
-1)+o(1)}\bigr) =
n^{1-\alpha+o(1)}.
\end{eqnarray*}
Also
\[
\mathbb E\bigl[\tilde{d}(v)^d \ind_{\bar{A}}\bigr] \leq
n^d \pr(\bar{A}) = O\bigl(n^{d-x}\bigr) = O
\bigl(n^{-\alpha}\bigr).
\]
Thus~(\ref{claimdvd}) holds, and hence so does~(\ref{eqn1ashow}), as
we noted earlier.

Now we may complete the proof using~(\ref{eqn1ashow}). We have
\begin{eqnarray*}
\mathbb EN_F(t_j,t_{j+1}) &=& \sum
_{k=1}^{\infty} \pr\bigl(\bigl
\{Z'_k \mbox{ fails}\bigr\} \cap\bigl
\{T'_k \le t_{j+1}\bigr\}\bigr)
\\
& \le& \sum_{k=1}^{\infty} \pr\bigl(\bigl
\{Z'_k \mbox{ fails}\bigr\} \cap\bigl
\{T'_k \le t_{j+1}\bigr\} \cap\{\tilde{T} \le
t_j\}\bigr)
\\
&&{} +  \sum_{k=1}^{\infty} \pr\bigl(\bigl
\{T'_k \le t_{j+1}\bigr\} \cap\{\tilde{T} >
t_j\}\bigr)
\\
& \le& (t_{j + 1} - t_j)n^{3 - d - \alpha+o(1)} +
N_j \pr(\tilde{T}>t_j)
\\
&&{}+ \mathbb E
\bigl[N_A(t_j,t_{j + 1})\ind_{N_A(t_j,t_{j + 1}) > N_j}
\bigr]
\\
& \le& n^{K+3-\alpha-d+o(1)} + o(1) = o(1).
\end{eqnarray*}
Thus $\mathbb E N_F(t_1,t_{3})=o(1)$, as required.


\subsection{Lower bounds: Showing \texorpdfstring{$\alpha$}{$alpha$} is $K$-bad}\label{lower2}

Here we want to prove that if we use the FDAR algorithm on a network
with $n$ nodes,
and the capacity $D \sim\alpha\ln n/\ln\ln n$ is not sufficiently
high, then many calls will be lost over an
interval of length $n^K$.
We shall use Lemma~\ref{lemblocking} inequality~(\ref{eqblock3}) to
obtain a lower bound on the
probability that a call $Z_k$ is lost. To use this lemma we need a
lower bound on
the number of saturated links at a vertex, and for this we need a lower
bound on the rate at which
calls arrive on a given link. Finally, to lower bound this rate we need
to upper bound the number
of saturated links, so that an arriving call wishing to use the link is
not too often blocked because
the ``partner'' link of the pair is saturated.
Thus to lower bound numbers of saturated links we must first upper
bound such numbers.

We say that a call $Z_k$ with endpoints $\{u,v\}$ and choices $j_1,j_2,
\ldots, j_d$
of intermediate nodes is \emph{blocked at} $u$ or \emph{blocked from} $v$
if the link $uj_1$
is saturated when the call arrives; that is, if $X_{T_k-}(uj_1) = D(n)$.
Clearly, if such a call $Z_k$ is not accepted onto a route, then in
particular, it is blocked at $u$ or $v$.
Also,
\[
\pr\bigl(\{Z_k \mbox{ blocked at } u\} \mid\mathcal
F_{T_k-},Z_k = \{u,v\}\bigr) = \frac{1}{n-2}\sum
_{j \neq u,v} \ind_{X_{T_k-}(uj) = D(n)}.
\]
Therefore, for each $v \in V$,
\begin{eqnarray*}
&& \pr\bigl(\{Z_k \mbox{ blocked from } v\} \mid\mathcal
F_{T_k-}, v \in Z_k\bigr)
\\
&&\qquad= \frac{1}{(n - 1)(n - 2)}\sum_{u \neq v}\sum
_{j
\neq u,v} \ind_{X_{T_k-}(uj) = D(n)}.
\end{eqnarray*}
Fix a node $v \in V$ and $0 < \delta< 1$. Then on the event
$S_{T_k-}(\mbox{at } u) \le(n-2) \delta/2$ for all nodes~$u$,
\[
\pr\bigl(\{Z_k \mbox{ blocked from } v\} \mid\mathcal
F_{T_k-}, v \in Z_k\bigr) \le\delta/2.
\]
In other words, while for each node $u$ the number of full links $uj$
is at most $(n-2) \delta/2$,
the probability that a new call which selects link $vj$ as its first
choice is blocked by the
``partner'' link $uj$ (where $u$ is the random other end of the call)
is at most $\delta/2$.
Thus, while for each node $u$ the number of full links $uj$ is at most
$(n-2) \delta/2$,
the arrival rate of calls onto each link $vj$ for $j \neq v$ is at
least $\lambda(1-\delta/2)$.

For\vspace*{1pt} $0 \leq s_0 \leq s_1$ let
$A'_{s_0,s_1}$ be the event that
$S_{t}(\mbox{at } u) \le(n-2)\delta/2$ for all nodes $u$ and all
times $t \in[s_0,s_1]$.
For each load vector $x$ and each node $j \neq v$, let $f_j(x)$ be the number
of calls in progress on the link $vj$. Also
let $\tilde{W}^{(vj)}$ be independent birth-and-death processes for $j
\neq v$,
each with arrival rate $\lambda_j = \lambda(1-\delta/2)$, death rate 1,
population 0 at time $0$, and population limit $n_j=D$.
Let $\tilde{W}_t (v) = \sum_j \ind_{\tilde{W}_t^{(vj)} = D}$,
the number of the $\tilde{W}^{(vj)}$ processes in state $D$ at time $t$.
Now we may apply Lemma~\ref{lemnew} on $[s_0,s_1]$, with $N=n-1$ and
$A$ as the event $A'_{s_0,s_1}$,
to obtain, for each integer $k \geq0$,
%
%
\begin{equation}
\label{eqnSDbound} \pr\bigl( \bigl\{ S_{s_1}(\mbox{at } v) \leq k \bigr
\}
\cap A'_{s_0,s_1} \bigr) \leq\pr\bigl(\tilde{W}_{s_1-s_0}
(v) \leq k \bigr).
\end{equation}
It follows that
\[
\pr\bigl( S_{s_1}(\mbox{at } v) \geq k \bigr) \geq\pr\bigl(
\tilde{W}_{s_1-s_0} (v) \geq k \bigr) - \pr\bigl(\overline{A'_{s_0,s_1}}
\bigr),
\]
and summing over $k=1,\ldots,n-1$ gives
\[
\mathbb ES_{s_1}(\mbox{at } v) \geq\mathbb E\tilde{W}_{s_1-s_0}
(v) - n \pr\bigl(\overline{A'_{s_0,s_1}}\bigr).
\]
It is well known that in equilibrium the $n-1$ immigration-death
processes $\tilde{W}^{(vj)}$ are i.i.d. random variables with a
Poisson distribution $\po(\lambda(1-\delta/2)) $ truncated
at~$D$.
Since, by standard theory, each $\tilde{W}^{(vj)}$
converges to equilibrium exponentially fast,
there exists a constant $\tilde{c} > 0$ such that, uniformly over
$t \geq\tilde{c} \ln n$ and $j \neq v$,
$\pr( \tilde{W}_t^{(vj)}= D) \ge n^{-\alpha+o(1)}$,
and so $\mathbb E[\tilde{W}_t (v)] \ge n^{1-\alpha+ o(1)}$.
Thus, assuming $s_1 \geq s_0+\tilde{c} \ln n$, for each vertex $v$, we have
%
%
\begin{equation}
\label{eqnESDbound} \mathbb ES_{s_1}(\mbox{at } v) \geq n^{1-\alpha+
o(1)} -
n \pr\bigl(\overline{A'_{s_0,s_1}}\bigr).
\end{equation}
Before we break into two cases as in the proof of the upper bound,
let us establish some more notation.
Let $t_0 = (K+5 + \tilde{c}) \ln n$, let $t_1 \ge t_0$, let $t'_1 =
t_1-\tilde{c} \ln n$,
and let $t_2=t_1 + n^K$. (For the lower bound proof, it is not
important to distinguish between the cases where $X_0$ is
stochastically at most the stationary process and where it is not.)
List the calls arriving after $t_1$ as $Z'_1, Z'_2, \ldots,$
arriving at times $t_1 < T'_1 < T'_2 < \cdots.$
As in the upper bound proof, $N_A(t_1,t_2)$ is the number of calls
arriving during the interval $(t_1,t_2]$,
and $N_F(t_1,t_2)$ is the number of calls that arrive during the
interval $(t_1,t_2]$
and are not accepted.

Suppose first that $0 < \alpha< \min\{1, (K+2)/d\}$.
(We consider the remaining case
$1 \leq\alpha< K+3-d$ later.)
Recall that $D \sim\alpha\ln n/\ln\ln n$.
Let $0 < \delta< \min\{ 1, (K+2)/d \} - \alpha$.
Using inequality~(\ref{binpo2})
\begin{eqnarray*}
\pr\bigl( \bigl\{ S_{t}(\mbox{at } v) \leq2n^{1-\alpha-\delta} \bigr\}
\cap A'_{s_0,s_1} \bigr) & \leq& \pr\bigl(
\tilde{W}_{t-s_0} (v) \le2n^{1-\alpha-\delta}\bigr)
\\
& \le& \exp\bigl(-n^{1-\alpha+ o(1)} \bigr)
\end{eqnarray*}
uniformly over nodes $v$ and times $t$ such that $s_0+ \tilde{c} \ln n
\leq t \leq s_1$.

For $0 \leq s_0 \leq s_1$ let $A_{s_0,s_1}$ denote the event that
$S_{t}(\mbox{at } v) \ge n^{1-\alpha- \delta}$ for all $v \in V$ and
all $t \in[s_0,s_1]$.
By the above and Lemma~\ref{lemstretch}(a), with $\tau= n^K$,
$a=b = n^{1-\alpha- \delta}$, and $\sigma= (2\lambda
)^{-1}n^{-\alpha- \delta}$,
\[
\pr\bigl(\overline{A_{t_1,t_2}} 
\cap A'_{t'_1,t_2}
\bigr) \leq\exp\bigl(-n^{1-\alpha- \delta+o(1)}\bigr) = o\bigl
(n^{-K-2}\bigr).
\]
Also, by Lemma~\ref{lemsat}, and Lemma~\ref{lemstretch}(a), with
$\tau=n^K + \tilde{c} \ln n$, $a=b = (n-2) \delta/4$,
and $\sigma= n^{-1/2}$,
\[
\pr\bigl( \overline{ A'_{t'_1,t_2}} \cap\bigl\{T \le
t'_1\bigr\} \bigr) = o\bigl(n^{-K-2}\bigr).
\]
Further by~(\ref{eqnTt2})
\[
\pr\bigl( T >t'_1\bigr) \leq\mathbb E\llVert
X_0 \rrVert_1 e^{-t'_1} = o\bigl(n^2
\ln n\bigr) \cdot e^{-(K+5)\ln n} =o\bigl(n^{-K-2}\bigr).
\]
It thus follows that
%
%
\begin{equation}
\label{eqnAbar-a} \pr\bigl(\overline{A'_{t'_1,t_2}}\bigr) = o
\bigl(n^{-K-2}\bigr)
\end{equation}
and
%
%
\begin{equation}
\label{eqnAbar} \pr(\overline{A_{t_1,t_2}}) = o\bigl(n^{-K-2}
\bigr).
\end{equation}

By Lemma~\ref{lemblocking}, equation~(\ref{eqblock3}), on the event
$B$ that
$S_{T'_k-}(\mbox{at } v) \ge(n-2)^{1-\alpha- \delta}$ for each $v
\in V$,
\[
\pr\bigl(Z'_k \mbox{ fails} \mid\mathcal
F_{T'_k-}\bigr) \ge\tfrac{1}2 (n-2)^{-(\alpha+ \delta)d}:=
p_0.
\]
Note that both $B$ and $\{T'_k \le t_2\}$ are in $\mathcal F_{T'_k-}$,
and so
\begin{eqnarray*}
\pr\bigl(Z'_{k} \mbox{ fails} \cap\bigl
\{T'_k \le t_2\bigr\}\bigr) & \geq& \pr
\bigl(Z'_{k} \mbox{ fails} \cap B \cap\bigl
\{T'_k \le t_2\bigr\}\bigr)
\\
& = & \mathbb E \bigl( \pr\bigl(Z'_k \mbox{ fails}
\mid\mathcal F_{T'_k-}\bigr) \ind_{B} \ind_{ \{T'_k \le t_2\}}
\bigr)
\\
& \ge& p_0 \pr\bigl(B \cap\bigl\{T'_k \le
t_2\bigr\}\bigr)
\\
& \ge& p_0 \pr\bigl(A_{t_1,t_2} \cap\bigl\{T'_k
\le t_2\bigr\}\bigr),
\end{eqnarray*}
where the last inequality follows since
\[
A_{t_1,t_2} \cap\bigl\{T'_k \le t_2
\bigr\} \subseteq B \cap\bigl\{T'_k \le t_2
\bigr\}.
\]
Now
\begin{eqnarray*}
\mathbb E\bigl[N_F(t_1,t_2)\bigr] & = &
\sum_{k=1}^{\infty} \pr\bigl(\bigl
\{Z_k' \mbox{ fails}\bigr\} \cap\bigl
\{T'_k \le t_2\bigr\}\bigr)
\\
& \ge& p_0 \sum_{k=1}^{\infty} \pr
\bigl(A_{t_1,t_2} \cap\bigl\{T'_k \le
t_2\bigr\}\bigr) = p_0 \mathbb E\bigl(
\ind_{A_{t_1,t_2}} N_A(t_1,t_2)\bigr).
\end{eqnarray*}
Let $N_0 = 2 \mathbb E[N_A(t_1,t_2)] = 2\lambda{n \choose2 } n^K$.
Note that by~(\ref{eqnAbar}) we have $N_0 \pr(\overline
{A_{t_1,t_2}}) = o(1)$, and, since $\mathbb E[N_A(t_1,t_2)] \le N_0/2$,
we have
$\mathbb E[N_A (t_1,t_2)\ind_{N_A(t_1,t_2) > N_0}] = o(1)$.
Thus
\begin{eqnarray*}
&&\mathbb E\bigl[\ind_{\overline{A_{t_1,t_2}}} N_A(t_1,t_2)
\bigr]
\\
&&\qquad =  \mathbb E\bigl[\ind_{\overline{A_{t_1,t_2}}} N_A(t_1,t_2)
\ind_{N_A(t_1,t_2) \le N_0}\bigr] + \mathbb E\bigl[\ind_{\overline
{A_{t_1,t_2}}}
N_A(t_1,t_2) \ind_{N_A(t_1,t_2) > N_0}\bigr]
\\
&&\qquad \le N_0 \pr(\overline{A_{t_1,t_2}}) + \mathbb E
\bigl[N_A(t_1,t_2) \ind_{N_A(t_1,t_2) > N_0}
\bigr] = o(1),
\end{eqnarray*}
and hence, for $n$ large enough,
\[
\mathbb E\bigl[N_F(t_1,t_2)\bigr] \ge
\tfrac{1}2 p_0 \mathbb E\bigl[N_A(t_1,t_2)
\bigr] \ge\tfrac{1}2 n^{K+2 - (\alpha+ \delta)d} = n^{\Omega(1)},
\]
as required.

Now consider the remaining case, when
$1 \le\alpha< K+3 -d$; see Figures~\ref{fig2} and~\ref{fig1}.
Recall that $\mathbb E[\tilde{W}_t (v)] \ge n^{1-\alpha+o(1)}$
uniformly over $t \ge\tilde{c} \ln n$ and $v \in V$.
By Lemma~\ref{lemblocking}, inequality~(\ref{eqblock3}),
\[
\pr\bigl(Z'_k \mbox{ fails} \cap\bigl
\{T'_k \le t_2\bigr\}\bigr) \ge
\frac{1}2 n^{-1} (n-2)^{-d} \sum
_v \mathbb E\bigl[ \bigl(S_{T'_k-}(\mbox{at } v)
\bigr)^d \ind_{\{T'_k \le t_2\}}\bigr].
\]
Since\vspace*{1.5pt} $S_{T'_k-}(\mbox{at } v) $ takes nonnegative values and we seek
a lower bound,
we may replace the exponent $d$ here by 1.
Let $S_t$ be the total number of saturated links at time $t$, so that
$S_t= \sum_v S_{t}(\mbox{at } v) $. Then by~(\ref{eqnESDbound})
and~(\ref{eqnAbar-a}),
for $t_1 < t \le t_2$,
\[
\mathbb ES_t \geq n \bigl(n^{1-\alpha+o(1)} - n \pr\bigl(\overline
{A'_{t_1',t_2}}\bigr)\bigr) = n^{2-\alpha+o(1)}.
\]
Hence, using the PASTA result Lemma~\ref{lempasta1} for the second
inequality below,
\begin{eqnarray*}
\mathbb E\bigl[N_F(t_1,t_2)\bigr] & = &
\sum_{k=1}^{\infty} \pr\bigl(Z'_k
\mbox{ fails} \cap\bigl\{T'_k \le t_2\bigr
\} \bigr)
\\
& \geq& \frac{1}2 n^{-1} (n-2)^{-d} \mathbb E\Biggl[
\sum_{k=1}^{\infty} \sum
_v S_{T'_k-}(\mbox{at } v) \ind_{ \{T'_k \le t_2\}}
\Biggr]
\\
& \geq& \frac{1}2 n^{-1} (n-2)^{-d} \mathbb E
\bigl[N_A(t_1,t_2)\bigr] \inf
_{t \in
(t_1,t_2]} \mathbb E[S_t]
\\
& \geq& n^{-1-d+2 +2-\alpha+o(1)} = n^{K+3-d-\alpha+o(1)} = n^{\Omega(1)},
\end{eqnarray*}
as required.

Now suppose $t_1 < t_0$ and $t_2 = t_1 + n^K$. Then we can apply the
above argument to calls arriving during the interval
$[t_0,t_2]$ with the same conclusion, and so $\mathbb E[N_F(t_1,t_2)] =
n^{\Omega(1)}$ in this case also, as required.


\section{Proof of Theorem~\texorpdfstring{\protect\ref{thmbdar}}{1.3}}\label{secproof2}

After introducing some notation and preliminary results, we will
discuss separately
the upper and lower bound parts of the theorem.

Fix an integer $d\ge2$ and a constant $K>0$.
Let $\phi= (101 + K)/\ln2$.
We choose times $t_0, t_1, t_2$, depending on $n$, as follows.
If $X_0$ is stochastically at most the equilibrium distribution, let
$t_0 \geq0$,
and otherwise let $t_0 \geq(K+8) \ln n$: now let
%
%
\begin{equation}
\label{eqndefntaub} t_1 = t_0+\phi\ln n\quad\mbox{and}\quad
t_2= t_1+ n^K.
\end{equation}
(Note that for convenience we have treated $t_0$ here slightly
differently from the statement of
Theorem~\ref{thmbdar}.)
Now fix a constant $0<\delta< 1$.
For each $t\in[t_0,t_2]$, let $A^0_t$ be the event
\[
\bigl\{(X_s(v)\le(1+ \delta)\lambda(n-1)\ \forall s
\in[t_0,t], \forall v \bigr\};
\]
by
Lemma~\ref{lemcalls3}, $\pr(\overline{A^0_{t_2}}) =O(n^{-K-3})$.

Also, let $A^1_t$ be the event
\[
\bigl\{ S_{s}^{}(\mbox{via } v) \le(n-2) \delta/4\
\forall s \in[t_0,t], \forall v \bigr\}.
\]
By Lemmas~\ref{lemsat}~and~\ref{lemstretch}
[with $a=b =(n-2)\delta/8$, $\tau= \phi\ln n + n^K$ and $\sigma= n^{-1/2}$],
\[
\pr\bigl(\overline{A^1_{t_2}} \cap\{\tilde{T} \le
t_0\}\bigr) =O\bigl(n^{-K-3}\bigr),
\]
and hence by~(\ref{eqnTt3}) also $\pr(\overline{A^1_{t_2}}) =O(n^{-K-3})$.

Recall that for each link $vw$, $X_t(vw)$ is the load of link $vw$
at time $t$, that is, the number of calls using this link at time $t$.
For $v \in V$ and $h=0,1, \ldots,$ let $L_t(v,h)$ be the number of
links $vw$ ($w \neq v$)
at $v$ with $X_t(vw) \ge h$ [so, in particular, for each $v \in V$,
$L_t(v, 0)=n-1$
for all $t$].
For $v \in V$ and $h=0,1, \ldots,$ we let $H_t(v,h)= \sum_{k \ge h}
L_t (v,k)$.

Let $c= \max\{c_1,c_2 \}$, where $c_1$ and $c_2$ are constants,
respectively, defined in Sections~\ref{secupper} and~\ref{seclower} below.
We will show that Theorem~\ref{thmbdar} holds with this value of $c$
and with $\kappa= K+7 + \phi$.


\subsection{Upper bound}\label{secupper}

Let the constant $c_1=c_1(\lambda, d, K)$ be as in~(\ref{capconst}) below,
and assume, as in the discussion preceding~(\ref{eqnTt3a}), that
${{\ln\ln n} \over{\ln d}}+c_1 \le D(n) = o(n)$, as $n \to\infty$.
We shall show that
a.a.s. no calls arriving during the interval $[t_1, t_2]$ of length $n^K$
fail.
We assume that $t_0 \geq(K+8) \ln n$ and $t_1= t_0+\phi\ln n$, as
at~(\ref{eqndefntaub}) above:
we will discuss briefly at the end of this subsection the case when
$X_0$ is stochastically
at most the equilibrium distribution, and we do not have a burn-in time
(so we allow then any $t_1 \geq0$).

Given a positive integer $h_0$, a decreasing
sequence of nonnegative numbers $(\alpha_h)_{h \ge h_0}$, and an increasing
sequence of times $(\tau_h)_{h\ge h_0}$ such that $t_0 \le\tau_h \le
t_1$ for each $h$, let
\[
B_t(h_0)= \bigl\{ L_s(v,h_0)
\le2 \alpha_{h_0}\ \forall s \in[\tau_{h_0},t],
\forall v \bigr\},
\]
and for $h=h_0+1,h_0+2, \ldots$ let
\[
B_t(h) = \bigl\{ H_s(v,h) \le2 \alpha_h\ \forall s \in[\tau_h,t], \forall v\bigr\}.
\]
Also, for each $h$, let $B(h) = B_{t_2}(h)$. Observe that if $B(h)$
holds, then
each link has load at most $h+2 \alpha_h -1$, at each time $t \in
[\tau_h,t_2]$.

The idea of the proof is to choose a sequence of
about $\ln\ln n/\ln d$ numbers $\alpha_h$ decreasing quickly from a
constant multiple of $n$ to zero,
and an increasing sequence of times $\tau_h$ for $h=h_0,h_0+1,h_0+2,
\ldots$ satisfying
$t_0 \le\tau_h \le t_1 $ for all~$h$. Then
the aim is to show that $B(h_0)$ holds a.a.s.,
and that, if $B(h)$ holds a.a.s., then so does $B(h+1)$, and to deduce that
$B(h)$ holds a.a.s. for some $h$ with $h+2\alpha_h \leq D$.
Thus a.a.s. no link is ever saturated during $[t_1,t_2]$, and so no call
can fail
during that interval.

We choose $h_0$ and a decreasing sequence of numbers $\alpha_h \geq0$
as follows. First, let
\[
h_0= \bigl\lceil\max\bigl\{ 8 \lambda, 768 \lambda^2
\bigr\} \bigr\rceil\quad\mbox{and}\quad \alpha_{h_0} = \min\biggl\{
\frac{n-1}{8},\frac{n-1}{768 \lambda}\biggr\}.
\]
Note that $\alpha_{h_0} \geq\lambda(n-1)/h_0$. Hence, on $A_{t_2}^0$,
for each $t \in[t_0, t_2]$,
since $X_t(v) \leq2 \lambda(n-1)$, we have
$L_t(v,h_0) \leq2 \lambda(n-1)/h_0 \leq2 \alpha_{h_0}$ and\vspace*{1pt} so
$A_{t_2}^0 \subseteq B(h_0)$.
We may (and do) assume that $n$ is sufficiently large that
$\alpha_{h_0} \geq14(K+4) \ln n$.
Next let the values $\alpha_h$ be defined by setting
%
%
\begin{equation}
\label{eqrecurrence} \frac{\alpha_h}{n-1} = 6\lambda\biggl( \frac
{8\alpha_{h-1}}{n-1}
\biggr)^d,
\end{equation}
for $h = h_0+1, h_0+2,\ldots, h^*$, where $h^*=h^*(n)$ is the largest
$h$ such that
$\alpha_h \ge14 (K+4) \ln n$.
(We shall see shortly that there is such an $h$.)
Also, define $\alpha_{h^*+1} = 14 (K+4)\ln n$ and $\alpha_{h^*+2} = 2K+7$.
Recurrence~(\ref{eqrecurrence}) can be rewritten as
%
%
\begin{equation}
\label{eqrecurrence-1} \tilde{\alpha}_h = 48\lambda\cdot\tilde{
\alpha}_{h-1}^d,
\end{equation}
where $\tilde{\alpha}_h = 8 \alpha_h/(n-1)$.
Since $\tilde{\alpha}_{h_0} \leq1$, it follows that for $h_0+1 \leq
h \leq h^*$
\[
\tilde{\alpha}_h = (48\lambda)^{1+d+\cdots+d^{h-h_0 -1}} \tilde{
\alpha}_{h_0}^{d^{h-h_0}} \leq(48\lambda\cdot\tilde{
\alpha}_{h_0})^{1+d+\cdots+d^{h-h_0 -1}}.
\]
But now, since $48\lambda\cdot\tilde{\alpha}_{h_0} \le\frac{1}2$,
for $h_0 \leq h \leq h^*$ we have
\[
\frac{8\alpha_h}{n-1}= \tilde{\alpha_h} \leq{ \biggl(\frac{1}2
\biggr)}^{(d^{h-h_0}-1)/(d-1)},
\]
and so $h^{*}(n) =\ln\ln n/\ln d + O(1)$. We now set
%
%
\begin{equation}
\label{capconst} c_1 = \sup_{k} \biggl\{ h^*(k) +
4K+16 - {{\ln\ln k} \over{\ln d}} \biggr\}
\end{equation}
so that
\[
D = D(n) \ge{{\ln\ln n} \over{\ln d}}+c_1 \geq h^*(n) +2+ 2 \alpha
_{h^* +2}.
\]
Now define an increasing sequence $(\tau_h)_{h \ge h_0}$ of times as follows.
Let
$\gamma_h = 48 \lceil\log_2 {(2\alpha_{h}/\alpha_{h+1})} \rceil$
for $h =h_0, \ldots, h^*-1$,
let $\gamma_{h^*} = 48 \log_2 n$,
and let $\gamma_{h^*+1} = (K+4) \log_2 n$.
Note that $2 \alpha_{h^*}/ \alpha_{h^* +1} \leq n/\ln n = o(n)$
and so $\gamma_{h^*} \ge\break 48 \lceil\log_2 {(2\alpha_{h^*}/\alpha
_{h^*+1})} \rceil$
for $n$ sufficiently large; this will be needed in the proof of
Lemma~\ref{lemupper1}.
Let $\tau_{h_0}=t_0$, and let $\tau_h=\tau_{h-1}+\gamma_{h-1}$ for
$h= h_0+1,h_0+2, \ldots, h^*+2$.
Thus $\tau_{h^*+2}= t_0 + \sum_{h=h_0}^{h^*+1} \gamma_h$.
Note that
\begin{eqnarray*}
\tau_{h^*}-t_0 & = & \sum_{h=h_0}^{h^*-1}
\gamma_h = 48 \sum_{h=h_0}^{h^*-1}
\bigl\lceil\log_2 {(2\alpha_h/\alpha_{h+1})}
\bigr\rceil
\\
& \le& 96\bigl(h^*-h_0\bigr)+48 \sum_{h=h_0}^{h^*-1}
(\log_2 \alpha_h -\log_2
\alpha_{h+1} )
\\
& \leq& 96 h^* + 48 \log_2{\alpha_{h_0}} \le49
\log_2 n
\end{eqnarray*}
for $n$ sufficiently large, and then
\[
\tau_{h^*+2}-t_0= \tau_{h^*} -t_0+
\gamma_{h^*} + \gamma_{h^*+1} \leq(101 + K) \log_2 n
= \phi\ln n. 
\]
Thus $\tau_{h^*+2} \leq t_1$.

As noted above, $A_{t_2}^0 \subseteq B(h_0)$, and so $\pr(B(h_0))$ is
near 1.
We shall show that $\pr(\overline{B(h)} \cap B(h-1))$ is small for
each $h=h_0+1,\ldots, h^{*}+2$,
which will yield that $\pr(B(h^*+2))$ is close to 1.
Hence, as we discussed earlier, since $D \geq h^*+2+2\alpha_{h^*+2}$,
a.a.s. throughout $[t_1,t_2 ]$, there are no full links. More precisely,
we shall show that
\[
\pr\bigl(\overline{B\bigl(h^*+2\bigr)}\bigr)=o\bigl(n^{-K-2}\bigr).
\]
Let $N_A(t_1,t_2)$ be the number of arrivals in $(t_1,t_2]$; then
$N_A(t_1,t_2) \sim\break \po(\lambda{n \choose2} (t_2-t_1))$. Let
$N_F(t_1,t_2)$ be the number of calls that fail during $(t_1,t_2]$.
Then
%
%
\begin{eqnarray}\label{ubshow}
\mathbb EN_F(t_1,t_2) & = & \mathbb
E[N_F \ind_{B(h^*+ 2)}] + \mathbb E[N_F
\ind_{\overline{B(h^*+ 2)}}]
\nonumber
\\
& \le& \lambda n^{K+2} \pr\bigl(\overline{B\bigl(h^*+2\bigr)}\bigr) +
\mathbb E\bigl[N_A (t_1,t_2)
\ind_{N_A(t_1,t_2) \ge\lambda n^{K+2}}\bigr]
\\
& = & o(1). \nonumber
\end{eqnarray}
This yields the desired upper bound of Theorem~\ref{thmbdar}
when the distribution of $X_0$ need not be stochastically dominated by
the stationary distribution.

To prove that $\pr(\overline{B(h)} \cap B(h-1))$ is small
for each $h$, we first show that if $B(h-1)$ holds, then a.a.s. for each
$v$ there exists a (random)
time $\tau_h(v) \in[\tau_{h-1},\tau_h]$ such that
$H_{\tau_h(v)}(v,h) \le\alpha_h$. We then show that a.a.s.
$H_t(v,h) \le2 \alpha_h$ for all $t \in[\tau_h (v),t_2]$ and all $v
\in V$.

For each node $v \in V$ and for each integer $h=h_0+1, \ldots, h^*+2$, let
\[
C(v,h)=\bigl\{ \exists\tau_h(v)\in[\tau_{h-1},
\tau_h]\dvtx H_{\tau
_h(v)}(v,h) \le\alpha_h \bigr\}.
\]
Let also $C(h) =\bigcap_v C(v,h)$, so that
\[
\overline{C(h)}=\bigl\{\exists w\dvtx  H_t(w,h)>
\alpha_h\ \forall t \in[\tau_{h-1},
\tau_h]\bigr\}
\]
is the event that there is a node $u$ such that the number of calls
with height at least $h$ at $u$ is greater than $\alpha_h$ throughout
$[\tau_{h-1},\tau_h]$.

%
\begin{lemma} \label{lemupper1}
\[
\sum_{h=h_0+1}^{h^*+2} \pr\bigl(\overline{C(h)}
\cap B(h-1)\bigr) = o\bigl(n^{-K-2}\bigr).
\]
\end{lemma}

\begin{pf}
The idea of the proof here is, for a fixed $v$ and $h$, to consider the
random value $H_t(v,h)$ at jump times $t$, when the value changes by
$0$ or $\pm1$.
We upper bound the probability of a positive change and lower bound the
probability
of a negative change, and then use Lemma~\ref{lemhittime1}.
For $h=h^* +2$ we need a slightly different argument, using Lemma~\ref{lembin}.

Fix a node $v$ and a height $h$ with $h_0+1 \le h \le h^*+1$.
Let $J_0(v)=\tau_{h-1}$, and enumerate the jump times of the process
of arrivals (possibly failing) and terminations of calls with one end $v$
after $J_0(v)$ as $J_1(v), J_2(v), \ldots.$
For $k=0,1,\ldots$ let $R_k = H_{J_k(v)}(v,h)$ and for
$k=1,2,\ldots$ let $Y_k= R_k- R_{k-1}$, so that
\[
R_k = R_0 + \sum_{j=1}^k
Y_j.
\]
Note that each jump $Y_k \in\{-1,0,1\}$ and is
$\mathcal F_{J_{k}(v)}$-measurable and hence also $\mathcal
F_{J_{k+1}(v)-}$-measurable,
and that the sum $\sum_{k\dvtx  \tau_{h-1} < J_k(v) \le\tau_h} Y_k$ is
the net change in $H_{t}(v,h)$ during the interval $(\tau_{h-1},\tau_h]$.
For $h=h_0, \ldots, h^*-1$, let
$m_h = \lceil12 \lambda n \rceil\lceil\log_2 (2\alpha_h/\alpha
_{h+1}) \rceil$,
which is ${\le}\frac{1}2 \gamma_h \lambda(n-1)$ for $n$ large enough.
Let also $m_{h^*} = \lceil12 \lambda n \rceil\lceil\log_2 n \rceil$,
which is ${\le}\frac{1}2 \gamma_{h^*} \lambda(n-1)$ for $n$ large enough.
Note that for each $h=h_0+1,\ldots,h^*+1$, we have
$J_{m_{h-1}}(v) \le\tau_h$ a.a.s., since by inequality~(\ref{binpo2}),
\[
\Pr\bigl(J_{m_{h-1}}(v) > \tau_h\bigr) \le\pr\bigl(\po\bigl(
\lambda(n-1) \gamma_{h-1}\bigr) < m_{h-1} \bigr) \le
e^{-\gamma_{h-1}\lambda(n-1)/8}.
\]
Now define events $E_k$ for $k=0,1, \ldots$ by letting
\[
E_k = A^0_{J_{k+1}(v)-} \cap B_{J_{k+1}(v)-}(h-1) =
A^0_{J_k(v)} \cap B_{J_k(v)}(h-1),
\]
and let $E=\bigcap_{k=0}^{m_{h-1}-1} E_k$.
We saw earlier that $\pr(\overline{A^0_{\tau_h}}) =O(n^{-K-3})$. Thus
\begin{eqnarray*}
\pr\bigl(\overline{E} \cap B(h - 1)\bigr) & \le& \pr\Biggl(\bigcup
_{k=0}^{m_{h - 1} - 1} \overline{A^0_{J_k(v)}}
\mbox{ for some } v \Biggr)
\\
& \le& \pr\bigl(J_{m_{h - 1}}(v) > \tau_h \mbox{ for some } v
\bigr)+ \pr\bigl(\overline{A^0_{\tau_h}}\bigr) = O
\bigl(n^{ -K - 3}\bigr).
\end{eqnarray*}

Now we obtain bounds for the probabilities (conditional on the past) of
jumps in $H_t (v,h)$:
an upper bound for the probability that a jump $Y_k$ is positive,
and a lower bound for the probability that a jump $Y_k$ is negative,
so that we can use Lemma~\ref{lemhittime1}.
On the event $E_{k-1}$,
upper bounding $\pr(J_k(v) \mbox{ is an arrival time} \mid\mathcal
F_{J_k(v)-})$ by~1,
we obtain for $n$ sufficiently large [since $n-2 \geq\frac{1}2(n-1)$
for $n \ge3$]
\begin{eqnarray*}
\pr(Y_k =1\mid\mathcal F_{J_k(v)-}) & \le& { \biggl(
{{2\max_w L_{J_k(v)-}(w,h-1)} \over{n-2}} \biggr)}^d
\\
& \le& { \biggl({{2\max_w H_{J_k(v)-}(w,h-1 )}\over{n-2}} \biggr)}^d
\\
& \le& { \biggl({{8 \alpha_{h-1} } \over{n-1}} \biggr)}^d \leq
\frac{\alpha
_h}{6 \lambda(n-1)},
\end{eqnarray*}
where the last inequality is from~(\ref{eqrecurrence}) and the choice
of $\alpha_{h^*+1}$.
Thus on the event~$E_{k-1}$,
\[
\pr(Y_k =1\mid\mathcal F_{J_k(v)-}) \leq a\qquad\mbox{where } a=
\frac
{\alpha_h}{6 \lambda(n-1)}.
\]

Now consider negative steps. The rate of arrivals of calls with one end
$v$ is $\lambda(n-1)$,
and on $A^0_{J_k(v)-}$ there are at most $2 \lambda(n-1)$ active calls
with one end $v$.
It follows that
on $A^0_{J_k(v)-}$,
\[
\pr(Y_k=-1\mid\mathcal F_{J_k(v)-}) \ge{{H_{J_k(v)-}(v,h)}\over{3
\lambda(n-1)}}=
{{R_{k-1}}\over{3\lambda(n-1)}},
\]
and so, for each $y \ge\alpha_h$, on $E_{k-1}\cap\{ R_{k-1}= y\}$,
\[
\pr(Y_k=-1\mid\mathcal F_{J_k(v)-})\ge{{y} \over{3\lambda(n-1)}}
= by\qquad\mbox{where } b = \frac{1}{3 \lambda(n-1)}.
\]
Note that for $t \ge\tau_{h-1}$, on $B_t(h-1)$, $H_t(w,h) \le2
\alpha_{h-1}-1$
for all nodes $w$. Note also that $\frac{2a}{b} = \alpha_h$.
Let $r = \alpha_h$, and let $\tilde{r}$ satisfy
$\alpha_h + 1 \le\tilde{r} \leq2 \alpha_{h-1}$.
Then by Lemma~\ref{lemhittime1}
\begin{eqnarray*}
&& \pr\bigl(E \cap\bigl\{ H_{J_k(v)}(v, h)>\alpha_h\
\forall k \le m_{h-1}\bigr\} \mid H_{J_0(v)}(v, h) = \tilde{r}
\bigr)
\\
&&\qquad\leq2e^{- \alpha_h /14} = O\bigl(n^{-K-4}\bigr).
\end{eqnarray*}
Summing over all nodes $v$, it follows that, uniformly over $h_0 +1 \le
h \le h^*+1$,
\[
\pr\bigl(\overline{C(h)} \cap B(h-1)\bigr) \le\pr\bigl(\overline{C(h)}
\cap E
\bigr) + \pr\bigl(\overline{E} \cap B(h-1)\bigr)= O\bigl(n^{-K-3}\bigr).
\]

Now let $h=h^*+2$.
We say that a call has height $h$ at $v$ if it is routed onto a link $vw$
that already has $h-1$ calls
at the time.
Let $N_H(v)$ denote the number of new calls for $v$ arriving during
$(\tau_{h^*+1},\tau_{h^*+2}]$
with height at least $h^*+2$. We will use Lemma~\ref{lembin} to show that
with high probability $N_H(v) \leq K+3$ for each
$v$.
Then we will see that with high probability no calls with height at
least $h^*+2$ at time
$\tau_{h^*+1}$ last until time $\tau_{h^*+2}$.

Enumerate calls with one end $v$ arriving after time $\tau_{h^*+1}$
as $Z'_1(v), Z'_2(v), \ldots$
with arrival times $J'_1(v), J'_2(v), \ldots.$
Recall that $\gamma_{h^*+1} = (K+4)\ln n$, and define
$m_{h^*+1} = \lceil2\gamma_{h^*+1}\lambda(n-1) \rceil$.
For $k=0,1, \ldots$ let
\[
E'_k= A^0_{J'_{k+1}(v)-} \cap
B_{J'_{k+1}(v)-}\bigl(h^*+1\bigr). 
\]
Further let $E'=\bigcap_{k=0}^{m_{h^*+1}-1} E'_k$.
For each $k=1,2, \ldots,$ on $E'_{k-1}$,
\[
\pr\bigl(Z'_k(v) \mbox{ has height } \ge h^*+2 \mid
\mathcal F_{J'_k(v)-} \bigr) \le p_1,
\]
where
\[
p_1= \biggl( \frac{4 \alpha_{h^*+1}}{n-2} \biggr)^d = \biggl(
\frac
{56 (K+4) \ln n}{n-2} \biggr)^d.
\]
Further we note that, for each positive integer $r$,
\[
\pr\bigl(B(m_{h^*+1},p_1) \ge r\bigr) \le(m_{h^*+1}p_1)^r
=O\bigl(\bigl(n^{-d+1} (\ln n)^{d+1}\bigr)^r
\bigr),
\]
where, as earlier, $B(n,p)$ is a binomial random variable with
parameters $n$ and $p$.

For $k=1, 2, \ldots,$ let $Y'_k$ denote \mbox{$\ind_{\{Z'_k(v)\ \mathrm{has}\
\mathrm{height}\,\ge\, h^*+2\}}$}.
Let $N_A(v)$ be the number of calls with one end $v$ arriving during the
interval $(\tau_{h^*+1}, \tau_{h^*+2}]$,
and let $N'_A(v)$ be the number of calls with one end $v$ arriving
during the
interval $(\tau_{h^*+1},t_2]$.
Then, using Lemma~\ref{lembin} (with $p=p_1$, $t=m_{h^*+1}$,
$Y_i=Y'_i$, $E_i=E'_i$,
$\mathcal F_i = \mathcal F_{J'_i(v)-}$ and $k=r$) for each integer $r
\ge K+4$,
\begin{eqnarray*}
&& \pr\bigl(\bigl\{N_H(v) \ge r\bigr\} \cap E'\bigr)
\\
&&\qquad  \le \pr\Biggl(\Biggl\{\sum_{k=1}^{m_{h^*+1}}
Y'_i \ge r\Biggr\} \cap E'\Biggr) + \pr
\bigl(N_A(v) > m_{h^*+1}\bigr)
\\
&&\qquad \le \pr\bigl(B(m_{h^* + 1},p_1) \ge r\bigr) + \pr\bigl(\po
\bigl(\lambda(n - 1)\gamma_{h^* +
1}\bigr) > m_{h^* + 1}\bigr) = o
\bigl(n^{ - K - 3}\bigr).
\end{eqnarray*}
Summing over all $v \in V$,
\begin{eqnarray*}
&& \pr\bigl(\bigl\{N_H(v) \ge r \mbox{ for some } v\bigr\} \cap B
\bigl(h^*+1\bigr)\bigr)
\\
&&\qquad \le \sum_v \pr\bigl(\bigl
\{N_H(v) \ge r\bigr\} \cap E'\bigr) + \pr\bigl(
\overline{E'} \cap B\bigl(h^*+1\bigr)\bigr)
\\
&&\qquad \leq o\bigl(n^{-K-2}\bigr) + \pr\bigl(\overline{A^0_{t_2}}
\bigr) + \sum_v \pr\bigl(N'_A(v)
< m_{h^*+1}\bigr) = o\bigl(n^{-K-2}\bigr).
\end{eqnarray*}
Also, on $B_{\tau_{h^*+1}}(h^*+1)$ there are at most $28 (K+4) \ln n$ calls
present at time $\tau_{h^*+1}$ with height at least $h^*+2$, and so
the probability
that at least one survives to time $\tau_{h^*+2}$ is at most
$28 (K+4) \ln n e^{-\gamma_{h^*+1}} = o(n^{-K-2})$. Hence
\[
\pr\bigl(\overline{C\bigl(h^*+2\bigr)} \cap B\bigl(h^*+1\bigr)\bigr)= o
\bigl(n^{-K-2}\bigr),
\]
as required.
\end{pf}

We now show that a.a.s., for each $h = h_0+1,\ldots, h^*+2$,
there will be no ``excursions'' that cross upwards from $\alpha_h$ to
at least
$2 \alpha_h$; that is, $H_t (v,h)$ cannot exceed $2 \alpha_h$ during
the time interval $(\tau_v (h), t_2]$ for any $v \in V$ and any
$h=h_0+1,\ldots, h^*+2$.

%
\begin{lemma} \label{lemupper2}
\[
\sum_{h=h_0+1}^{h^*+2} \pr\bigl(\overline{B (h)}
\cap B(h-1) \cap C(h)\bigr) = o\bigl(n^{-K-2}\bigr).
\]
\end{lemma}

\begin{pf}
Take $h \in\{h_0, \ldots, h^*+2\}$.
The only possible start times for an upward crossing excursion of $H_t
(v,h)$ are
arrival times during $[\tau_{h-1},t_2]$.
Let $N_0 = 2 \lambda n^{K+1}$.
Then the probability that, for some $v$, more than $N_0$ calls with one
end $v$ arrive
during the interval $(\tau_h(v),t_2]$ is $O(n^{-K-3})$.

Now consider a fixed node $v$.
Let $J_0 = \tau_h(v)$, and let $J_1, J_2,\ldots$ be the jump times of
the process of
arrivals (possibly failing) and terminations of calls with one end $v$
after time $\tau_h(v)$.
For $k=0,1, \ldots$ let $R_k = H_{J_k}(v,h)$,
and for $k=1,2, \ldots$ let $Y_k = R_k-R_{k-1}$.
Then each $Y_k \in\{-1,0,1\}$ and is
$\mathcal F_{J_k}$-measurable and thus also $\mathcal
F_{J_{k+1}-}$-measurable. For $k=0,1, \ldots,$ let
\[
E_k = A^0_{J_{k}} \cap B_{J_{k}} (h-1).
\]
As in the proof of Lemma~\ref{lemupper1}, on $E_{k-1}$
\begin{eqnarray*}
\pr(Y_k =1 \mid\mathcal F_{J_k}) \le q_h^+&:=&
\biggl( \frac{8 \alpha
_{h-1}}{n-1} \biggr)^d,
\\
\pr(Y_k =-1 \mid\mathcal F_{J_k}) \ge q_h^-&:=& \frac{\alpha_h}{3
\lambda(n-1)},
\end{eqnarray*}
and for $h \le h^*+1$,
\[
q_h^+ \leq{{\alpha_h} \over{6 \lambda(n-1)}}=q_h^-/2.
\]

Let $p=q_h^+$, $q=q_h^-$ and $a=\lfloor\alpha_h \rfloor-1 \geq
\alpha_h -2$.
By Lemma~\ref{lemcrossing}, the probability that the event $A^0_{t_2}
\cap B(h-1)$
occurs and any given excursion during $(\tau_h(v),t_2]$ leads
to a ``crossing'' is at most $(q_h^+/q_h^-)^{\alpha_h -2}$,
and so for $h=h_0+1,\ldots, h^*+2$, summing over all $v \in V$ and over
all possible excursion starting times,
\begin{eqnarray*}
&& \pr\bigl(\overline{B(h)} \cap B(h-1) \cap C(h) \bigr)
\\
&&\qquad \le 2\lambda n^{K+2} \bigl(q_h^+/q_h^-
\bigr)^{\alpha_h -2}+ \pr\bigl(\bar{A}^0_{t_2}\bigr)
+ O\bigl(n^{-K-3}\bigr).
\end{eqnarray*}
For $h=h_0+1,\ldots, h^*+1$,
\[
\bigl(q_h^+/q_h^-\bigr)^{\alpha_h -2}
\leq2^{-\alpha_h +2} = 2^{-14 (K+4) \ln n +2}= O\bigl(n^{-2K-5}\bigr),
\]
and so the above bound is $O(n^{-K-3})$.
For $h=h^*+2$,
\[
\frac{q_h^+}{q_h^-} = \biggl( \frac{8 \cdot14 (K+4) \ln
n}{n-1} \biggr)^d \cdot
\frac{\lambda(n-1)}{2K+7} =O\bigl(n^{1-d} \ln^{d} n\bigr) = O
\bigl(n^{-1} \ln^2n\bigr)
\]
and so
\[
\bigl(q_h^+/q_h^-\bigr)^{\alpha_h-2} = O
\bigl(n^{-2K-5}\bigr) \cdot\ln^{O(1)}n,
\]
and the lemma follows.
\end{pf}

We may now complete the proof of Theorem~\ref{thmbdar}(a).
Recall that $A^0_{\tau_2} \subseteq B_{h_0}$, and that
$\pr(\overline{A^0_{\tau_2}})=o(n^{-K-2})$. Then
\begin{eqnarray*}
&& \pr\bigl(\overline{B\bigl(h^*+2\bigr)}\bigr)
\\
&&\qquad \leq \pr\bigl(\overline{A^0_{\tau_2}}\bigr)+ \pr\bigl(
\overline{B(h_0)} \cap A^0_{\tau
_2}\bigr) + \sum
_{h=h_0+1}^{h^*+2} \pr\bigl(\overline{B(h)} \cap
B(h-1)\bigr)
\\
&&\qquad = \pr\bigl(\overline{A^0_{\tau_2}}\bigr) + \sum
_{h=h_0 + 1}^{h^* + 2} \pr\bigl(\overline{C(h)} \cap B(h - 1)
\bigr)
\\
&&\quad\qquad{} + \sum_{h=h_0 + 1}^{h^* + 2} \pr\bigl(
\overline{B(h)} \cap C(h) \cap B(h - 1)\bigr)
\\
&&\qquad =  o\bigl(n^{-K-2}\bigr).
\end{eqnarray*}

This completes the proof of~(\ref{ubshow}) and thus of the upper bound
of Theorem~\ref{thmbdar},
for the case of general starting configuration.

Finally let us consider the case when the distribution of the initial
state $X_0$
is stochastically dominated by the stationary distribution $\pi$.
Let us set $t_0=0$ and consider $t_1 \in[0,(k+8)\ln n)$.\vspace*{1pt}

Let $B_t = \bigcap_{h=h_0}^{h^*+2} B'_t (h)$, where the events
$B'_t(h)$ are like the events $B_t(h)$ above,
but with $2$ replaced by $3/2$. We may adapt the upper bound proof
described above to show that
$B_t$ holds a.a.s. for $t$ large enough.
Thus $B_0$ must hold a.a.s. for the equilibrium distribution. But $B_0$ is
a decreasing event, and so
$B_0$ must hold a.a.s. for any initial distribution stochastically at most
the equilibrium distribution.
Now we deduce as above that, for all $v$ and all
$t \in[0, t_1 + n^K]$, $L_t(v,h_0) \leq2 \alpha_{h_0}$, and
$H_t(v,h) \leq2 \alpha_h$ for each $h=h_0+1,\ldots,h^*+2$.
Finally, we may deduce as before that the expected number of calls that
fail during
$[0, t_1 + n^K]$ is $o(1)$, and this completes the proof.


\subsection{Lower bound}\label{seclower}

Let the constant $c_2=c_2(\lambda, d, K)$ be as defined
below, and let $D = D(n) \le\frac{\ln\ln n}{\ln d} -c_2$. Let
$0 < \varepsilon< \min\{1, (K+2)/d \}$.
Once again, we work on the interval $[t_1,t_2]$ of length $n^K$ defined
in~(\ref{eqndefntaub}).
We shall show that
a.a.s. for each $v$ at least $(n-1)^{1-\varepsilon}$ links $vw$ incident
on $v$ are saturated (and so unavailable) throughout the interval,
and hence
a.a.s. at least
$ n^{K +2-\varepsilon d-o(1)}$ calls arriving during the interval fail.

Given a sequence of nonnegative numbers $(\alpha_h)_{h \ge0}$ and a
sequence of times $(\tau_h)_{h \ge0}$ such that $t_0 \le\tau_h \le
t_1$ for each $h$, let
\[
B_t(h) = \bigl\{L_s (v,h) \ge\alpha\ \forall s
\in[\tau_h,t], \forall v\bigr\},
\]
and let $B(h) = B_{t_2}(h)$.
We shall choose numbers
$\alpha_0,\alpha_1,\ldots,$ starting with $\alpha_0=n-1$ and
decreasing rapidly.
We shall further choose an increasing sequence of
times $\tau_h$, $h = 0,1,\ldots,$ such that
$t_0 \le\tau_h \le t_1$ for each $h$.
Our aim is to show that $B (D(n))$ occurs a.a.s., with
a value $\alpha_{D(n)} \ge(n-1)^{1-\varepsilon}$,
so that there are always many saturated links.

The\vspace*{1pt} numbers $\alpha_h$ are given as follows.
Let $\nu= \frac{\min\{1,\lambda\}}{24 e^d}$, so that $0 < \nu<1$.
Now let $\alpha_{0} = n-1$, and for $h=1,2,\ldots$ define $\alpha_h$
by setting
%
%
\begin{equation}
\label{eqRecurrence} \frac{\alpha_h}{n-1} = \frac{\nu}{h} \biggl(
\frac{\alpha_{h-1}}{n-1} \biggr)^d.
\end{equation}
Since\vspace*{1pt} $\frac{1}{12} \le e-1$, it is easily checked that $2 \alpha_h
\le\alpha_{h-1} (1-e^{-1})$,
and so $(\alpha_h - 2\alpha_{h+1})^d \ge(\alpha_h/e)^d$ for each $h$.

We want to choose the constant $c_2$ in the upper bound on $D(n)$ above
such that for $n$ sufficiently large
\[
\alpha_{D(n)} \ge(n-1)^{1-\varepsilon}.
\]
To see that such a choice is possible, let $\beta_h = \frac{\alpha_h}{n-1}$.
Then $\beta_0 =1$ and
%
%
\begin{equation}
\label{eqRecurrence-1} \beta_h = \frac{\nu}{h} \beta_{h-1}^d
\qquad\mbox{for }h=1,2,\ldots.
\end{equation}
It follows that for each positive integer $h$,
%
%
\begin{equation}
\label{eqrecurrence-sol} \beta_h = \frac{\nu^{1+d+ \cdots+ d^{h-1}}}{
\prod_{i=1}^h i^{d^{h-i}}}.
\end{equation}
To upper bound the denominator in~(\ref{eqrecurrence-sol}), note that
for some $c_3 > 0$,
\[
\ln\bigl(h (h-1)^d (h-2)^{d^2}\cdots2^{d^{h-2}}
\bigr) = d^h\sum_{i=2}^{h}
d^{-i} \ln i \le c_3 d^h
\]
and so $\prod_{i=1}^h i^{d^{h-i}} \le e^{c_3 d^h}$.
It follows that for each $h \in\nats$,
\[
\beta_h \ge e^{-d^h (\ln(1/\nu) + c_3)}.
\]
Let $c_4$ be such that $d^{-c_4}( \ln(\frac{1}{\nu}) + c_3) \le
\varepsilon$;
if $h \le\ln\ln(n-1) /\ln d -c_4$, then
\begin{eqnarray*}
\beta_h & \ge& \exp\bigl(-\bigl(\ln(n-1)\bigr) d^{-c_4}
\bigl( \ln(1/\nu) + c_3\bigr) \bigr)
\\
& \ge& \exp\bigl(-\varepsilon\ln(n-1)\bigr) = (n-1)^{-\varepsilon}.
\end{eqnarray*}
Since $\ln\ln n \le\ln\ln(n-1) +1$ for $n$ large enough, we can
take $c_2 = c_4+1$.

For\vspace*{2pt} $h = 0,1,\ldots$ let
$\gamma_h = \frac{4}{\max\{1,\lambda\} (h+1)}$.
Now define an increasing sequence of times $\tau_h$ as follows.
Let $\tau_0=t_0$, and for $h=1,\ldots,$ let
$\tau_h=\tau_{h-1} +\gamma_{h-1}$.
Then
\[
\tau_{D(n)} - t_0 = \sum_{h=0}^{D(n)-1}
\gamma_h \le4 \sum_{h=1}^{D(n)}
\frac{1}h \le4 \ln(D+1)= O(\ln\ln\ln n).
\]
It follows that $\tau_D \le t_1$ for $n$ sufficiently large.

Since $\alpha_0=n-1$, it follows that $\pr(B(0))=1$; we prove
by induction that $ \pr(\overline{B(h)}) = O (n^{-K-3})$
for $h=1,\ldots, D(n)$, so that a.a.s. throughout $[t_1,t_2]$ for
each $v$ there are at least $(n-1)^{1-\varepsilon}$ saturated links $vw$
incident on $v$. The main step is to show that $\pr(\overline{B(h)}
\cap B(h-1))$ is small for each $h$; to do this,
we first show that if $B(h-1)$ occurs, then a.a.s. for each $v$ there
exists a time
$\tau_h(v) \in[\tau_{h-1},\tau_h]$ such that $L_{\tau_h(v)}(v,h)
\geq2
\alpha_h$.

For each node $v \in V$ and each positive integer $h$, let
\[
C(v,h) = \bigl\{ \exists\tau_h(v) \in[\tau_{h-1},\tau
_h]\dvtx L_{t_h(v)}(v,h) \ge2 \alpha_h \bigr\},
\]
and let $C(h)=\bigcap_v C(v,h)$.

\begin{lemma}
\label{lemlower1}
%
%
\begin{equation}
\label{eqnCbarB} \sum_{h=1}^D \pr\bigl(
\overline{C(h)} \cap B(h-1)\bigr) = o\bigl(n^{-K-2}\bigr).
\end{equation}
\end{lemma}

\begin{pf}
The idea is very similar to that in the proof of Lemma~\ref{lemupper1}. We consider the variable $L_t(v,h)$ at jump times $t$,
when it changes by $0$ or $\pm$1. We lower bound the probability of a
positive change and upper bound the probability of a negative change,
and use a reversed version of Lemma~\ref{lemhittime}.

Fix a node $v$ and an integer $h \ge1$.
Let $J_0 (v) = \tau_{h-1}$ and enumerate the jump times of the process
of arrivals (possibly failing) and terminations of calls with one end
$v$ after time $J_0 (v)$ as
$J_1(v), J_2 (v), \ldots.$ For $k=0,1,\ldots$ let
$R_k = L_{J_k(v)}(v,h)$ and for $k=1,2, \ldots$ let
$Y_k= R_k - R_{k-1}$, so that
\[
R_k = R_0 + \sum_{j=1}^k
Y_j.
\]
Then each $Y_k \in\{-1,0,1\}$, is $\mathcal F_{J_k(v)}$ and hence also
$\mathcal F_{J_{k+1}(v)-}$-measurable,
and $\sum_{k\dvtx  \tau_{h-1} < J_k (v) \le\tau_h} Y_k$ is the net
change in
$L_t (v,h)$ during $(\tau_{h-1},\tau_h]$.

For $h=0,1, \ldots,$ let $m_h = 2\min\{1,\lambda\}(n-1)/(h+1)= \frac
{1}2 \lambda(n-1) \gamma_h$.
Note that,
for $h =0,1, \ldots, D(n)$,
\[
\pr\bigl(J_{m_{h-1}} (v) > \tau_h \bigr) \le\pr\bigl(\po
\bigl(\lambda(n-1) \gamma_{h-1}\bigr) < m_{h-1}\bigr) \le
e^{-\gamma_{h-1}
\lambda(n-1)/8}.
\]
For $k=0,1, \ldots$ let
$E_k = A^0_{J_{k+1}(v)-} \cap A^1_{J_{k+1}(v)-} \cap B_{J_{k+1}(v)-}(h-1)$.
Let $E = \bigcap_{k=0}^{m_{h-1}-1}E_k$.
Recalling that $\pr(\overline{A^0_{t_2}\cup A^1_{t_2}}) = O(n^{-K-3})$,
\begin{eqnarray*}
\pr\bigl(\overline{E} \cap B(h-1)\bigr) &\le&\pr\bigl(J_{m_{h-1}}(v) >
\tau_h \mbox{ for some } v\bigr) + \pr\bigl(\overline{A^0_{\tau_h}
\cup A^1_{\tau_h}}\bigr)
\\
&=& O\bigl(n^{-K-3}\bigr).
\end{eqnarray*}
We now seek a lower bound (conditional on the past) on the probability
that the jump $Y_k$ takes value 1.
First note that on $A^0_{J_k(v)-}$ the conditional probability that
$J_k(v)$ is an arrival time (for a call for $v$) is at least $1/(2
+\delta) \ge\frac{1}3$. Now note that $Y_k$ takes value 1 if the $k$th
call (with endpoints $v$ and $u$, for any random choice of $u \neq v$)
is routed onto a link $vw$ with load exactly $h-1$ at $J_{k-1}(v)$;
this will happen if, in particular, for every intermediate node $w_i$
selected, the link $vw_i$ has load exactly $h-1$ at time $J_{k-1}(v)$,
and at least one of the
``partner'' links $uw_i$ is not blocked at time $J_{k-1}(v)$.

Now we want to consider $u$ picked uniformly at random (u.a.r.) from $V
\setminus\{v\}$
and $w_1,\ldots,w_d$ picked u.a.r. from $V \setminus\{v,u\}$. We may
pick $u$ and $w_1,\ldots,w_d$ as follows.
First pick $w_1$ u.a.r. from $V \setminus\{v\}$, then pick $u$ u.a.r. from
$V \setminus\{v,w_1\}$,
then pick $w_2,\ldots,w_d$ independently and u.a.r. from $V \setminus\{
v,u\}$.
[This\vspace*{1.5pt} gives exactly the same distribution on the $(d+1)$-tuple
$u,w_1,\ldots,w_d$.]
On\vspace*{1.5pt} $A^1_{J_k(v)-}$ we have $S_{{J_k(v)-} }^{}(\mbox{via } w) \le(n-2)/2$
for all nodes $w$; and so, whatever $w_1$ is picked,
the probability conditional on $\mathcal F_{J_k(v)-}$ that $uw_1$ is
saturated is at most $\frac{1}2$.
Hence, on~$A^0_{J_k(v)-} \cap A^1_{J_k(v)-}$,
\begin{eqnarray*}
&& \pr(Y_k=1 \mid\mathcal F_{J_k(v) -})
\\
&&\qquad \geq \frac{1}3 \frac{L_{J_k(v) -} (v,h - 1) - L_{J_k(v) -}
(v,h)}{n-1}
\\
&&\qquad\quad{}\times \frac{1}2 { \biggl(
\frac{L_{J_k(v) -} (v,h - 1) - 1 - L_{J_k(v) -} (v,h)}{n-2} \biggr)}^{d
- 1}
\\
&&\qquad \geq \frac{1}6 { \biggl( {{L_{J_k(v)-}(v,h-1) - 1 -
L_{J_k(v)-}(v,h)} \over{n
- 1}}
\biggr)}^d.
\end{eqnarray*}

It follows that, on $E_{k-1} \cap(R_{k-1} < 2 \alpha_h)$,
\[
\pr(Y_k = 1\mid\mathcal F_{J_k(v) - }) \ge\frac{1}6
\biggl(\frac{ \alpha_{h -
1} - 2 \alpha_h }
{n - 1} \biggr)^d \ge\frac{e^{-d}}{6} \biggl(\frac{ \alpha_{h - 1}} {n - 1} \biggr)^d
= {{4h\alpha_h } \over{\min\{1, \lambda\} (n - 1)}}
\]
by~(\ref{eqRecurrence}).
Thus on the event $E_{k-1} \cap(R_{k-1} < 2 \alpha_h)$,
\[
\pr(Y_k=1\mid\mathcal F_{J_k(v)-}) \ge2p\qquad\mbox{where } p
= {{2h\alpha_h } \over{\min\{1,\lambda\} (n-1)}}.
\]
Now we consider negative steps.
The probability that $J_k (v)$ is a
departure time of a given
call with one end $v$ is at most
$
\frac{1}{\lambda(n-1)}$,
and so
\[
\pr(Y_k=-1 \mid\mathcal F_{J_k(v)-}) \leq
{{h(L_{J_k(v)-}(v,h)-L_{J_k(v)-}(v,h+1))}\over{\lambda(n-1)}}.
\]
It follows that, for each $y < 2 \alpha_h$, on $E_{k-1} \cap(R_{k-1}=y)$,
\begin{eqnarray*}
\pr(Y_k=-1\mid\mathcal F_{J_k(v)-}) & \le&
{{hy} \over{\lambda(n-1)}}\le{{2h\alpha_h } \over{\lambda(n-1)}}
\\
& \leq& {{2h\alpha_h } \over{\min\{1,\lambda\} (n-1)}}=p.
\end{eqnarray*}
Let $r_1 = 2\alpha_h$, and let $r_0$ be any positive
integer less than $2 \alpha_h$. Note that $q_h^+ m_{h-1} \ge4
\alpha_h \ge2 (r_1 - r_0)$. By a natural ``reversed'' version of
Lemma~\ref{lemhittime}
(i.e., by Lemma~7.2 in~\cite{lmc03a}), for any value of $r_0 \le r_1$,
\begin{eqnarray*}
&& \pr\bigl(E \cap\bigl(L_{J_k(v)} (v,h) < 2\alpha_h\ \forall k
\in\{1,\ldots,m_{h - 1}\}\bigr) \mid L_{J_0(v)}(v,h -
1)=r_0\bigr)
\\
&&\qquad\leq e^{-\alpha_h /7} \leq e^{-\alpha_D /7} \le
e^{-\Omega(n^{1-\varepsilon})}.
\end{eqnarray*}

Note that we used Lemma~\ref{lemhittime1} in place of Lemma~\ref
{lemhittime} in the corresponding part of the proof of the upper
bound in Theorem~\ref{thmbdar}. The reason for this is that, in the
upper bound, we had to bring the quantity $H_t(v,h)$ from at most $2
\alpha_{h-1}$ to $\alpha_h$ (rather than from at least 0 to $2 \alpha
_h$) and, for large $h$, $\alpha_{h-1}$ and $\alpha_h$ are of a
different order of magnitude in $n$, and we did not want to
``give away'' the extra downward drift of $H_t(v,h)$ in the vicinity of
$\alpha_{h-1}$.

Summing over all $v$ we see that
\[
\pr\bigl(\overline{C(h)} \cap B(h-1)\bigr)\le\pr\bigl(\overline{C(h)}
\cap E
\bigr) + \pr\bigl(\overline{E} \cap B(h-1)\bigr) = O \bigl(n^{-K-3}\bigr).
\]
Thus we have now completed the proof of~(\ref{eqnCbarB}).
\end{pf}

We now need to prove that for each $h=1,2, \ldots, D(n)$, a.a.s. there
will be no
excursions that cross downwards from $2 \alpha_h$ to less than $\alpha_h$;
that is, each of the numbers $L_t (v,h)$ is unlikely to drop below
$\alpha_h$ during $(\tau_v(h),t_2]$.

\begin{lemma}
\label{lemlower2}
\[
\sum_{h=1}^{D} \pr\bigl(\overline{B(h)}
\cap B(h-1) \cap C(h)\bigr) = o\bigl(n^{-K-2}\bigr).
\]
\end{lemma}

\begin{pf}
Take $h \in\{1, \ldots, D(n)\}$ and $v \in V$. The only possible
start times for a crossing of $L_t(v,h)$ from $2\alpha_h$ to $\alpha_h$
are completion times of calls with one end $v$ during $[\tau
_{h-1},t_2]$. Let $N_0 = 4 \lambda n^K$.
Then
the probability that, for some $v$, more than $N_0$ calls with one end
$v$ terminate during the interval $(\tau_h(v),t_2]$ is $O(n^{-K-3})$.

Now consider a fixed node $v$.
Let $J_0 = \tau_h(v)$, and let $J_1, J_2, \ldots$ be the jump times
of the process of arrivals
(possibly failing) and completions of calls with one end $v$ after time
$\tau_h(v)$.
For $k=0,1, \ldots,$ let $R_k = L_{J_k}(v,h)$ and for $k=1,2, \ldots$
let $Y_k = R_k-R_{k-1}$.
Then each $Y_k \in\{-1,0,1\}$ and is $\mathcal F_{J_k}$ and hence
$\mathcal F_{J_{k+1}-}$-measurable.
For $k=0,1, \ldots$ let
\[
E_k= A^0_{J_{k+1}-} \cap A^1_{J_{k+1}-}
\cap B_{J_{k+1}-}(h-1).
\]
As in the proof of Lemma~\ref{lemlower1}, on $E_{k-1}$,
\begin{eqnarray*}
\pr(Y_k=1\mid\mathcal F_{J_k(v)-}) &\ge& q_h^+:
= {{4h\alpha_h } \over
{\min\{1,\lambda\} (n-1)}},
\\
\pr(Y_k=-1\mid\mathcal F_{J_k(v)-}) & \le&
q_h^-:= {{2h\alpha_h } \over{\min\{1,\lambda\} (n-1)}}
\end{eqnarray*}
and $q_h^-= \frac{1}{2} {q_h^+}$.

Analogously to the proof of Lemma~\ref{lemupper2}, we may
apply a reversed version of Lemma~\ref{lemcrossing} with
$p =q_h^-$, $q = q_h^+$, $a=\lfloor\alpha_h \rfloor-1$.
The probability that
the event $A^0_{t_2} \cap A^1_{t_2} \cap B(h-1)$ occurs and any given
excursion during
$(\tau_h(v), t_2]$ leads to a ``crossing'' is at most
$(q_h^-/q_h^+)^{\lfloor\alpha_h \rfloor-1} \le(1/2)^{\alpha_h-2}$.
Summing over all $v \in V$, for $h=1, \ldots, D(n)$,
\begin{eqnarray*}
\pr\bigl(\overline{B(h)} \cap B(h - 1) \cap C(h) \bigr) & \le& n N_0
\bigl(\tfrac{1}2\bigr)^{\alpha_h-2} + O\bigl(n^{- K - 3}\bigr)
+ \pr\bigl(\overline{A^0_{t_2}} \cup\overline{A^1_{t_2}}
\bigr)
\\
& \le& 4\lambda n^{K+2}\bigl(\tfrac{1}2\bigr)^{\alpha_D-2}
+ O\bigl(n^{ - K - 3}\bigr)
\\
&=& O\bigl(n^{ - K
- 3}\bigr).
\end{eqnarray*}\upqed
\end{pf}

Now, as in the proof of the upper bound,
\begin{eqnarray*}
\pr\bigl(\overline{B\bigl(D(n)\bigr)}\bigr) & \leq& \pr\bigl(\overline{B(0)}
\bigr) + \sum_{h=1}^{D(n)} \pr\bigl(
\overline{B(h)} \cap B(h - 1)\bigr)
\\
& = & \sum_{h=1}^{D(n)} \pr\bigl(
\overline{C(h)} \cap B(h - 1)\bigr) + \sum_{h=1}^{D(n)}
\pr\bigl(\overline{B(h)} \cap C(h) \cap B(h - 1)\bigr)
\\
& = & o\bigl(n^{-K-2}\bigr).
\end{eqnarray*}

As before, let $N_A(t_1,t_2)$ be the total number of calls arriving in
$(t_1,t_2]$; then
$N_A(t_1,t_2) \sim\po(\lambda{n \choose2} (t_2-t_1))$.
Also, as before, $N_F(t_1,t_2)$ is the number of calls that are lost
during $(t_1,t_2]$.
On the event $B_{T'_k-} (D(n))$, for $n$ sufficiently large,
\[
\pr\bigl(Z'_k \mbox{ fails} \mid\mathcal
F_{T'_k-}\bigr) \ge\frac{1}2 \biggl(\frac{(n-1)^{1-\varepsilon} -1}{n-2}
\biggr)^d \geq\frac{1}4 n^{-\varepsilon d}:=
p_1.
\]
Let $N_1 = \lceil\frac{1}2 \lambda{n \choose2} n^K \rceil$.
Let
$b^* = \frac{1}{32} \lambda n^{K+2-d \varepsilon}$, and let $B^*= \{
N_F(t_1,t_2) < b^*\}$.
Then, by Lemma~\ref{lembin},
\[
\pr\bigl(B^*\bigr) \le\pr\bigl(\overline{B\bigl(D(n)\bigr)}\bigr) + \pr
\bigl(N_A(t_1,t_2) < N_1\bigr)
+ \pr\bigl(B (N_1,p_1) < b^*\bigr)= o
\bigl(n^{ - K - 2}\bigr).
\]
Now suppose $0 \le t_1 \le t_0$, and let $t_2 = t_1 + n^K$. Then we can
apply the above argument to $[t_0,t_2]$ with the same conclusion.
Since $\varepsilon$ can be chosen arbitrarily small, this completes the
proof of the lower bound of Theorem~\ref{thmbdar}.


\section{Concluding remarks}\label{secconc}

We have considered the performance of two algorithms for a
continuous-time network routing problem,
strengthening and extending the earlier results in~\cite{lu99}
and~\cite{aku}, with full proofs.

For simplicity we have assumed throughout that the underlying network
is a complete graph,
but our results carry over in a straightforward way to a suitably
``dense'' subnetwork.
Consider, for example, the upper bound in Theorem~\ref{thmbdar} part~(a).
Let $\delta>0$, and suppose that, in the network with $n$ nodes, for
each pair of nodes $u$ and $v$
the number of possible intermediate nodes is at least $\delta n$.
[For instance, if] $0<p<1$ is fixed and the $n(n-1)$ possible links appear
independently with probability $p$,
then with high probability each pair of distinct nodes has
$(p^2+o(1))n$ common neighbours.]
Minor alterations to the proof of Theorem~\ref{thmbdar} part~(a) show
that we obtain the same conclusion:
if $D(n) \ge{\ln\ln n}/ {\ln d} +c$, and we use the BDAR algorithm,
then the expected number of failing calls during an interval of length
$n^K$ is $o(1)$. The only difference is that now the constant $c$
depends also on $\delta$.
Note that the leading term $\ln\ln n/ {\ln d}$ depends only on the
problem size $n$ and the number $d$ of choices, and not on $\delta$
(or on $\lambda$ or $K$).

For the dense networks we have been considering, it has been natural to
work with two-link routes.
If we wish to consider routing in sparser networks, for example, a
random graph as above but with $p=o(1)$,
then it would be natural to consider longer routes for calls, but we do
not pursue that here.


The analysis in~\cite{lu99} (see also~\cite{l00}) suggested that the
performance of the model could be upper
and lower bounded by differential equations. While that analysis was
nonrigorous, it turns out that a suitable
differential equation approximation, and concentration of measure
bounds, can indeed be obtained:
the details appear in~\cite{l08+}.
The main challenge was to disentangle the complex dependencies within
subsets of links to obtain a tractable
asymptotic approximation for the generator of the underlying Markov process.


\section*{Acknowledgements}

The authors are grateful for helpful comments from Benjamin Stemper and
a careful referee.




%

\printaddresses
\end{document}